\documentclass[10pt]{amsart}

\usepackage{amsmath}
\usepackage{amssymb}
\usepackage{latexsym}

\newcommand{\Zint}{\mathbb {Z}}    
\newcommand{\Rat}{\mathbb {Q}}     
\newcommand{\Rea}{\mathbb {R}}      
\newcommand{\Cplx}{\mathbb {C}}     
     
\newcommand{\Ha}{\mathfrak {h}}    
\newcommand{\halmos}{\rule{5pt}{5pt}}
\newcommand{\rank}{\mathrm {rank}\,}
\newcommand{\dist}{\mathrm {dist}}

\numberwithin{equation}{section}
\newtheorem{df}{\bf Definition}
\newtheorem{prop}{\bf Proposition}[section]
\newtheorem{thm}[prop]{\bf Theorem}
\newtheorem{lemma}[prop]{\bf Lemma}
\newtheorem{suble}[prop]{\bf Sublemma}

\newenvironment{rmk}{\noindent{\bf Remark}\hskip 5pt}{\hfill{$\Box$}}

\begin{document}
\baselineskip 21pt

\title[perturbation of Calogero-Moser-Sutherland system]
{The perturbation of the quantum Calogero-Moser-Sutherland system and related results}
\author{Yasushi Komori}
\address{Institute of Physics, University of Tokyo, Komaba, Tokyo 153-8902, Japan.}
\email{komori@gokutan.c.u-tokyo.ac.jp}
\thanks{Y.K. is supported by JSPS Research Fellowship for Young Scientists} 
\author{Kouichi Takemura}
\address{Department of Mathematical Sciences, Yokohama City University, 22-2 Seto, Kanazawa-ku, Yokohama 236-0027, Japan.}
\email{takemura@yokohama-cu.ac.jp}


\begin{abstract}
The Hamiltonian of the trigonometric Calogero-Sutherland model coincides 
with a certain limit of the Hamiltonian of the elliptic Calogero-Moser model.
In other words the elliptic Hamiltonian is a perturbed operator of the trigonometric one.

In this article we show the essential self-adjointness of the Hamiltonian of
the elliptic Calogero-Moser model and the regularity (convergence) of
the perturbation for the arbitrary root system.

We also show the holomorphy of the joint eigenfunctions of the
commuting Hamiltonians w.r.t the variables $(x_1, \dots ,x_N)$ for the
$A_{N-1}$-case.  As a result, the algebraic calculation of the
perturbation is justified.
\end{abstract}

\maketitle

\section{Introduction}

The Hamiltonian of elliptic Calogero-Moser model is given as follows (\cite{OP}),
\begin{equation}
\label{Hamiltonian_A}
\mathcal{H} :=- \frac{1}{2} \sum_{i=1}^{N} 
\frac{\partial ^{2}}{\partial x_{i}^{2}}
+ 
\beta (\beta -1)
\sum_{1 \leq i<j \leq N}
\wp ( x_{i}-x_{j}).
\end{equation}
where $\beta $ is the coupling constant.

This Hamiltonian reduces to that of the trigonometric Calogero-Sutherland model by setting $\tau \rightarrow \sqrt{-1} \infty$, where $\tau$ is the ratio of two basic periods of the elliptic function.

As for the trigonometric Calogero-Sutherland model, 
it is well-known for specialists that 
their eigenstates are given by the Jack polynomials (or the $A_{N-1}$-Jacobi polynomials).
So far, many researchers have studied the Jack polynomials 
and its $q$-deformed version, the Macdonald polynomials,
and clarified various properties such as
 the orthogonality, the norms, the Pieri formula, the Cauchy formula, 
and
the evaluations at $(1, \dots, 1)$.
The Calogero-Sutherland model is extended to 
those associated with simple Lie algebras.
From this point of view
the Hamiltonian \eqref{Hamiltonian_A} 
is called the $A_{N-1}$-type.
Studies of these models are being developed 
by using their algebraic structures. 

In contrast with the trigonometric models,
the elliptic models are less investigated and
the spectrum or the eigenfunctions are not sufficiently analyzed.
There are, however, some important progress due to Felder and Varchenko.
They clarify that the Bethe Ansatz works well
for the $A_{N-1}$-type elliptic Calogero-Moser model (\cite{FVthr}). 
Although this method may have applications to the spectral problem 
and indeed some partial results are obtained in the article (\cite{R,Tak}),
we will employ another approach in the present article.

In this article, we will add some knowledges of 
the elliptic Calogero-Moser model based 
on the analysis of the trigonometric model, 
which we will explain below.

One topic is the essential self-adjointness 
of the elliptic Calogero-Moser model 
for the arbitrary root system. 
Firstly we will establish it for the trigonometric model 
by taking the Jacobi polynomials as its domain in 
the space of square integrable functions.
We will obtain the elliptic version 
by perturbation.

Second topic is 
to obtain the eigenvalues and the eigenfunctions 
of the elliptic Calogero-Moser model for arbitrary root systems.
There are at least two ways to perform it.
The one is to use the Bethe Ansatz method, 
which is valid only for the $A_{N-1}$ case. 
From this viewpoint some results are obtained in (\cite{R,Tak}).
The other is to use the well-developed perturbation theory, 
which we will consider in this article. 
We regard the Hamiltonian of the elliptic Calogero-Moser model 
as the perturbed operator of the Calogero-Sutherland model 
by the parameter $p=\exp(2\pi\tau \sqrt{-1})$.
We have so abundant knowledge about 
the eigenvalues and the eigenfunctions 
of the Calogero-Sutherland model that
we can apply the perturbation method.
Then we will obtain the eigenvalues and the eigenfunctions 
as a formal power series of $p$.
In general, such formal power series does not converge. 
For example, consider the operator $H:=-\frac{d^2}{dx^2}+x^2+\alpha x^4$, 
then the formal power series of the eigenvalues and eigenfunctions 
diverge for any $\alpha \neq 0$.
However, in our cases, the formal power series converges 
if $p$ is sufficiently small.
The convergence is assured by the functional analytic 
method introduced by Kato and Rellich.
We mean the convergence of the eigenfunctions 
in the $L^2$-norm sense.

The other topics are the holomorphy of the eigenfunctions, 
the relationship with the higher commuting operators, 
and giving the elliptic analogue of the Jacobi polynomials, 
which are valid for the $A_{N-1}$ case.
The Kato-Rellich method does not give the holomorphy a priori.
We will obtain the holomorphy by using several properties 
of the Jack (or the $A_{N-1}$-Jacobi) polynomials.
Thanks to the holomorphy, the eigenspaces of the second-order Hamiltonian are 
compatible with the higher commuting operators.
By considering the joint eigenfunctions of the Hamiltonian and 
the higher commuting operators, 
we see the well-definedness of an elliptic analogue 
of the $A_{N-1}$-Jacobi polynomials.

We remark that Langmann obtained the algorithm for constructing 
the eigenfunctions and eigenvalues as the formal power series of $p$. 
(\cite{Lang})
His algorithm would be closely related to our one which is explained 
in section \ref{secpert}.

There are some merits for the perturbation method comparing 
to the Bethe Ansatz method. 
The calculation of the perturbation does not essentially 
depend on the coupling constant $\beta $
though the calculation of the Bethe Ansatz method 
strongly depends on $\beta $.
In addition,
the Bethe Ansatz method is applied to the $A_{N-1}$ type 
and $\beta \in \Zint_{>1}$ cases, 
but the perturbation method may be valid for all types and the coupling constant does not need to be an integer.

\section{Jacobi polynomials and self-adjointness}

The Hamiltonians of the trigonometric and the elliptic 
Calogero-Sutherland models are respectively given by
\begin{align}
\label{eq:Hamiltonian_trig} 
\mathcal{H}_{{T}}
&:=
-\Delta
+\sum_{\alpha\in{R}_+}k_\alpha(k_\alpha-1)|\alpha|^2
\Bigr(\frac{1}{4\sin^2(\langle\alpha,h\rangle/2)}-\frac{1}{12}\Bigr),
\\
\label{eq:Hamiltonian_ell} 
\mathcal{H}_{{E}}
&:=
-\Delta
+\sum_{\alpha\in{R}_+}k_\alpha(k_\alpha-1)|\alpha|^2
\wp(\langle\alpha,h\rangle),
\end{align}
where the coupling constant $k_\alpha$ is real and
invariant under the action of the Weyl group $k_\alpha=k_{w\alpha}$,
and $\wp(x)=\wp(x;\pi,\pi \tau)$ is the Weierstrass $\wp$ function.
For our later convenience,
we have subtracted a constant term from the original trigonometric
Hamiltonian.
$\Delta$ is the Laplacian on 
$T:=\mathfrak{h}_{\mathbb{R}}/2\pi{Q}{}^\vee$.
By using the variable $p=\exp(2\tau\pi\sqrt{-1})$,
we often write $\mathcal{H}_E(p)=\mathcal{H}_E$ in order to emphasize 
the dependency of $p$. 
In this notation, we have $\mathcal{H}_T=\mathcal{H}_E(0)$. 

We first show that the Hamiltonian of the trigonometric model
 is defined on a dense subspace of
$L^2(T,d\mu)^{W}$
where $\mu$ is the normalized Haar measure,
and is essentially self-adjoint with respect to
the inner product
\begin{equation}
  (f,g):=\int_T \overline{f}\cdot g.
\end{equation}
We denote by $\|\cdot\|:=(\cdot,\cdot)^{1/2}$
the norm in $L^2(T,d\mu)$.
We define $\mathcal{H}_{T}$
and $\mathcal{H}_{E}$ on $C^2(T)^W\cap \mathsf{D}(V)$, where $\mathsf{D}(V)$ denotes
the domain of the multiplication operators in
\eqref{eq:Hamiltonian_trig} and
\eqref{eq:Hamiltonian_ell}.
Then we see that these operators are symmetric.

If $k_{\alpha}\geq 2$,
$\mathcal{H}_{T}$ has a $C^2$-class 
$W$-invariant eigenfunction
\begin{equation}
  \mathcal{H}_{T}\Delta=E_0 \Delta,
\end{equation}
where 
\begin{align}
E_0&=(\rho(k)|\rho(k))-e_0, \\
e_0&=\frac{1}{12}\sum_{\alpha \in R_+} k_{\alpha }(k_{\alpha }-1)|\alpha|^2, \\
\rho(k)&=\frac{1}{2}\sum_{\alpha\in R_+}k_\alpha \alpha, \\
\Delta&=\prod_{\alpha\in{R}_+}
|\sin(\langle\alpha,h\rangle/2)|^{k_\alpha}.
\label{gstate}
\end{align}

Let $\mathbb{C}[P]$ be the polynomial ring of the weight lattice $P$.
For each $\lambda$, let $e^\lambda$ denote the corresponding element, 
so that $e^\lambda e^\mu=e^{\lambda+\mu}$ and $e^{0}=1$.
We also regard the element $e^\lambda$ as 
a function on $T$ by the rule 
$e^{\lambda}(\Dot{h}):=e^{\sqrt{-1}\langle\lambda,h\rangle}$
where $\Dot{h}\in T$ is the image of 
$h\in \mathfrak{h}_{\mathbb{R}}$.
Let $m_\lambda$ for $\lambda\in P_+$ be the monomial symmetric functions
\begin{equation}
  m_\lambda:
=|W_\lambda|^{-1}\sum_{w\in W}e^{w\lambda}
=\sum_{\mu\in W\lambda}e^\mu,
\end{equation}
where $W_\lambda$ denotes the stabilizer of $\lambda$ in $W$.
The set $\{m_\lambda|\lambda\in P_+\}$ forms a basis of $\mathbb{C}[P]^W$.
Define the partial order $\prec$ in $P$ by
\begin{equation}
  \nu\preceq \mu \Leftrightarrow \mu-\nu\in Q_+.
\end{equation}
Let $\|\cdot\|_\Delta$ and $(\cdot,\cdot)_\Delta$
denote the norm and the inner product in $L^2(T,\Delta^2 d\mu)$
respectively. 
\begin{df}[Heckman-Opdam]
\label{def:HO}
There exists a family of polynomials $\{J_\mu|\mu\in P_+\}$
which consists of a basis of $\mathbb{C}[P]^W$ satisfying
the following conditions:
\begin{gather}
J_\mu=m_\mu+\sum_{\nu\prec\mu} u_{\mu\nu}m_\nu, \\
(J_\mu,J_\nu)_\Delta=0,\qquad\text{if }\mu\neq\nu.
\end{gather}
\end{df}
Let $H_0:=\Delta^{-1} \mathcal{H}_{T} \Delta$.
Then these polynomials are characterized by the operator $H_0$
as its eigenfunctions.
\begin{prop}
\label{prop:eigen_H_0}
\begin{equation}
H_0 J_\mu=E_\mu J_\mu,
\end{equation}
where $E_\mu=(\mu+\rho(k)|\mu+\rho(k))-e_0$.
\end{prop}

It is well known that
the normalized Jacobi polynomials $\Tilde{J}_{\lambda }$ $(\lambda \in P_+)$ 
form a complete orthonormal system 
in the space
$L^2(T,\Delta ^2 d\mu)^{W}$ 
with the inner product $(\cdot,\cdot)_{\Delta}$
if $k_\alpha\geq 0$.
It follows that 
\begin{lemma} \label{conb}
Assume $k_{\alpha}\geq 0$. Then
$\mathcal{P}:=\mathbb{C}[P]^{W}\Delta$
is a dense subspace in
$L^2(T,d\mu)^{W}$. 
\end{lemma}

\begin{thm}
\label{thm:Hamiltonian_trig_is_SA}
Assume $k_{\alpha}\geq 2$. Then
$\mathcal{H}_{T}$ is essentially self-adjoint on 
$\mathcal{P}$.
\end{thm}
\begin{proof}
From Proposition \ref{prop:eigen_H_0} 
and $\mathcal{P}\subset C^2(T)^W\cap \mathsf{D}(V)$
we see that
$J_\lambda \Delta$ are the eigenfunctions of $\mathcal{H}_{T}$.
Then the theorem is obtained from Lemma \ref{conb}
since it implies that the range of $(\mathcal{H}_{T}\pm i)$ is dense.
\end{proof}
If $0<k_\alpha<2$, then $\Delta\not\in C^2(T)^W\cup\mathsf{D}(V)$
and $\mathcal{P}$ is not an appropriate domain for $\mathcal{H}_T$.
However Theorem \ref{thm:Hamiltonian_trig_is_SA} is generalized in the following sense
in terms of the adjoint operator $\mathcal{H}^*_T$:
\begin{thm}
Assume $k_{\alpha}\geq 0$. 
Then $\mathcal{H}^*_{T}|_\mathcal{P}$ is essentially self-adjoint on 
$\mathcal{P}$.
\end{thm}

We rewrite
$\mathcal{H}_{E}(p)=\mathcal{W}(p)+\mathcal{H}_{T}$, where
$\mathcal{W}(p)=(\mathcal{H}_{E}(p)-\mathcal{H}_{T})$
is a multiplication operator with 
\begin{equation}
\mathcal{W}(p)(h):=
  \sum_{\alpha\in{R}_+}k_\alpha(k_\alpha-1)|\alpha|^2
  \Bigl(
  \wp(\langle\alpha,h\rangle)-\frac{1}{4\sin^2(\langle\alpha,h\rangle/2)}
  +\frac{1}{12}
  \Bigr).
\end{equation}
By the formula (\ref{wpth}), we see that
\begin{equation}
\|\mathcal{W}(p)u\|\leq \|\mathcal{W}(p)\|_{\max}\|u\|,
\end{equation}
since the function $\mathcal{W}(p)(h)$
is a continuous function on $T$.
This implies that $\mathcal{W}(p)$ is 
bounded.
Hence we have
$\mathcal{H}^*_{E}(p)=\mathcal{W}(p)^*+\mathcal{H}^*_{T}$.
\begin{thm}[\cite{Kom}]
\label{thm:Hamiltonian_ell_is_SA}
Let $-1<p<1$ and $k_{\alpha}\geq 0$.
Then 
$\mathcal{H}^*_{E}(p)|_\mathcal{P}$ 
is essentially self-adjoint on $\mathcal{P}$.
\end{thm}
\begin{proof}
The symmetry of the operator $\mathcal{W}(p)$ is trivial.
Then 
we deduce
that 
$\mathcal{W}(p)+\mathcal{H}^*_{T}|_\mathcal{P}$ 
is essentially self-adjoint 
on $\mathcal{P}$. 
\end{proof}

In the next section, 
we abuse
the symbols $\mathcal{H}_T$ and $\mathcal{H}_E(p)$ 
for $\mathcal{H}^*_T|_\mathcal{P}$ and
$\mathcal{H}^*_E(p)|_\mathcal{P}$.


\section{Perturbation in the $L^2$-space} \label{sec:pertL2}

In this section, 
we employ the variable $p$ with $|p|<1$
instead of $\tau$ as a parameter of perturbation
and treat mainly the gauge-transformed Hamiltonian defined below
with $k_\alpha> 0$.
For a linear operator $T$,
we denote by $\mathsf{D}(T)$ its domain 
and by $\mathsf{R}(T)$ its range respectively.

\subsection{The resolvent in the $L^2$ space} \label{l2}
For a bounded linear operator $A$, we denote by $\| A \|$ the operator norm, 
\emph{i.e.}, $\| A \|:= \sup_{\| v  \|_{\Delta}=1 }\| Av\|_{\Delta}$.
We set
\begin{align}
W(p)&:= \Delta^{-1} (\mathcal{H}_E(p)-\mathcal{H}_T) \Delta(=\mathcal{W}(p)), \label{eq:Wp} \\
T(p)&:= \Delta^{-1} \mathcal{H}_E(p) \Delta.
\end{align}
Then $T(p)=H_0 + W(p)$  
is a closable operator on
$L^2(T,\Delta ^2 d\mu)^{W}$
with $\mathsf{D}(T(p))=\mathbb{C}[P]$, 
and particularly
if $-1<p<1$,
$T(p)$ is
an essentially self-adjoint operator.
Here
$W(p)$ is a bounded operator on
$L^2(T,\Delta ^2 d\mu)^{W}$ with an upper bound,
\begin{equation}
\|W(p)\|\leq W_{\max}(p):=4\sum_{n=1}^{\infty}\frac{n|p|^n}{1-|p|^n} \cdot \sum_{\alpha \in R_+} k_{\alpha }(k_{\alpha }-1)|\alpha |^2,
\end{equation}
which is monotonous with respect to $|p|$ and
 tends to 0 as $p\to 0$.


Let $\Tilde{T}$ denote the closure of a closable operator $T$.
Then $\Tilde{T}(p)$ for $p\in (-1,1)$ is
the unique extension of $T(p)$ to the self-adjoint operator.
In particular $\Tilde{H}_0=\Tilde{T}(0)$ 
is the self-adjoint extension of $H_0$.
$\Tilde{T}(p)$ is a self-adjoint holomorphic family \cite{Kat}.
Notice that the spectrum of the operator $\Tilde{H}_0$ is discrete.
Let $\sigma (\Tilde{H}_0)$ 
be the set of the spectrum and 
let $\rho(\Tilde{H}_0)$ be the resolvent set
of the operator $\Tilde{H}_0$.
We have 
\begin{equation} 
\sigma (\Tilde{H}_0)= \{(\lambda +\rho (k) | \lambda +\rho (k))-e_0| \lambda \in P_+ \}.
\label{CSev}
\end{equation}

The following proposition is obvious.
\begin{prop} \label{prop:multH0}
For each $a \in  \sigma (\Tilde{H}_0)$, the corresponding eigenspace $\{ v \in L^2(T , \Delta ^2 d\mu )^{W}$ $|~ \Tilde{H}_0v=a v\}$ is finite dimensional.
\end{prop}
For $\zeta \in \rho(\Tilde{H}_0)$,
the resolvent $(\Tilde{H}_0- \zeta )^{-1}$ is compact 
and 
$\| (\Tilde{H}_0- \zeta )^{-1}\|=(\dist(\zeta, \sigma (\Tilde{H}_0)))^{-1}$.
We have
\begin{equation}
 (\Tilde{H}_0- \zeta )^{-1}\sum_{\lambda }c_{\lambda} J_{\lambda}= \sum_{\lambda }(E_{\lambda }- \zeta )^{-1}c_{\lambda} J_{\lambda},
\end{equation}
where $\sum_{\lambda }c_{\lambda} J_{\lambda} \in  L^2(T,\Delta ^2 d\mu)^{W}$ and $\Tilde{H}_0J_{\lambda}=E_{\lambda }J_{\lambda}$.
The proof of the Kato-Rellich theorem also implies
the compactness of the resolvent of $\Tilde{T}(p)$ for $-1 <p<1$.

If $\|(\Tilde{H}_0- \zeta )^{-1} W(p) \|<1$, 
then $(\Tilde{H}_0-\zeta)^{-1}(\Tilde{T}(p)-\zeta)=1+(\Tilde{H}_0-\zeta)^{-1}W(p)$
has a bounded inverse 
by Neumann series and thus $\Tilde{T}(p)-\zeta$ has also a bounded inverse
\begin{align}
(\Tilde{T}(p)-\zeta )^{-1}
& = \bigl(1+(\Tilde{H}_0-\zeta )^{-1}W(p)\bigr)^{-1}(\Tilde{H}_0-\zeta )^{-1} \label{Neus} \\
& = \sum_{j=0}^{\infty}\bigl(-(\Tilde{H}_0-\zeta )^{-1}W(p)\bigr)^j (\Tilde{H}_0-\zeta )^{-1}. \notag
\end{align}
In particular,
 if $\| (\Tilde{H}_0- \zeta )^{-1} \| < \|W(p)\|^{-1}$, 
then the bounded inverse of $\Tilde{T}(p)-\zeta$ exists.
The right hand side of this expression implies
that the resolvent of $\Tilde{T}(p)$ is also compact.
By the equality 
$\|(\Tilde{H}_0- \zeta )^{-1}\|=(\dist(\zeta, \sigma (\Tilde{H}_0)))^{-1}$ 
and the equation \eqref{CSev}, 
the resolvent set $\rho(\Tilde{T}(p))$ is included outside 
the union of the closed disks 
$\dist(\zeta, \sigma (\Tilde{H}_0))\leq \|W(p)\|$.


\begin{prop} \label{proj}
Let $T$ be a closed operator 
with the resolvent set $\rho(T)$. 
Let $\Gamma_1 , \Gamma_2$ be circles 
which are contained in $\rho(T)$ and whose interiors are disjoint.
We set 
$$ P_i:=-\frac{1}{2\pi \sqrt{-1}}\int_{\Gamma_i}(T-\zeta )^{-1} d\zeta,\; \; (i=1,2).
$$
Then we have 
$P_i^2=P_i$ and $P_1P_2=P_2P_1=0$.
\end{prop}

\begin{prop} \label{proj2}
Let $P$, $Q$ be bounded operators 
subject to 
$P^2=P$, $Q^2=Q$ and $\| P-Q \|<1$. 
Then we have $\rank P= \rank Q$.
\end{prop}
\begin{proof}
For $u\in\mathsf{R}(Q)$, we have $Pu=u+Pu-Qu$ due to $Qu=u$.
Then
\begin{equation}
\|Pu\|\geq(1-\|P-Q\|)\|u\|,  
\end{equation}
which implies that 
$P|_{\mathsf{R}(Q)}:\mathsf{R}(Q)\to\mathsf{R}(PQ)\subset\mathsf{R}(P)$ 
is one-to-one and $\rank P\geq\rank Q$.
Similarly we have $\rank P\leq\rank Q$ and hence $\rank P= \rank Q$.
\end{proof}

Let $\Gamma\subset \rho(\Tilde{H}_0)$ be a circle
and 
let $r=\dist(\Gamma,\sigma(\Tilde{H}_0))$.
Then there exists $p_0>0$ such that for all $|p|<p_0$,
$\|W(p)\|<r$ and thus $\Gamma\subset \rho(\Tilde{T}(p))$.
Notice $\|(\Tilde{H}_0-\zeta)^{-1}\|\leq r^{-1}$ on $\Gamma$.
Let
$$ P_\Gamma(p):=-\frac{1}{2\pi \sqrt{-1}}\int_\Gamma(\Tilde{T}(p)-\zeta )^{-1} d\zeta.
$$
Then we have
\begin{align}
\| P_\Gamma(p)-P_\Gamma(0) \|
&\leq \frac{1}{2\pi}
\int_{\Gamma}\|(\Tilde{T}(p)-\zeta )^{-1}-(\Tilde{H}_0-\zeta)^{-1}\||d\zeta|
\\
&\leq \frac{1}{2\pi}
\int_{\Gamma} \sum_{j=1}^{\infty}
\|(\Tilde{H}_0-\zeta)^{-1}\|^{j+1}\|W(p)\|^j|d\zeta|
\\
&<
\Bigl(
\frac{1}{2\pi} 
\int_{\Gamma} |d\zeta|
\Bigr)
\frac{r^{-2}\|W(p)\|}{1-r^{-1}\|W(p)\|}.
\end{align}
Fix $a_i \in \sigma (\Tilde{H}_0)$.
Since the set $\sigma (\Tilde{H}_0)$ is discrete, 
we can choose a circle $\Gamma _i$ and $0<p_i$ 
such that $\Gamma _i$ contains only one element 
$a_i$ inside it 
and $P_i(p)=P_{\Gamma_i}(p)$ satisfying
$\| P_i(p)-P_i(0) \|<1$ for $|p|<p_i$. 
By Propositions \ref{proj}, \ref{proj2}, we 
see that $\rank P_i(p)= \rank P_i(0)$ and 
in particular, $P_i(p)$ is a degenerate operator.
By the proof of Proposition \ref{proj2},
we see that $V_i(p):=\mathsf{R}(P_i(p))$ is spanned by the image of 
$V_i(0)=\mathsf{R}(P_i(0))$, \emph{i.e.}, the eigenspace of $\Tilde{H}_0$ 
with the eigenvalue $a_i$;
We choose a basis of $V_i(p)$ as 
$
\{ 
P_i(p) \Tilde{J}_{\lambda}~|~\lambda\text{ such that }
\Tilde{H}_0\Tilde{J}_{\lambda}=a_i\Tilde{J}_{\lambda}\}
$, where $\Tilde{J}_{\lambda}$ is the normalized Jacobi polynomial.
One sees that $V_i(p)$ 
is a finite dimensional invariant subspace of 
$\Tilde{T}(p)$ due to the commutativity of
$P_i(p)$ and $\Tilde{T}(p)$.
\begin{lemma}
The matrix elements of $\Tilde{T}(p)|_{V_i(p)}:V_i(p)\to V_i(p)$
with respect to $P_i(p)\Tilde{J}_\lambda$
are real-holomorphic functions of $p$.
\end{lemma}
\begin{proof}
We define the functions $c^\mu_\lambda(p)$ 
and $d^\mu_\lambda(p)$ 
by
\begin{align}
\Tilde{T}(p)P_i(p)\Tilde{J}_\lambda
&=\sum_\mu c^\mu_\lambda(p)P_i(p)\Tilde{J}_\mu,\\
P_i(0)P_i(p)\Tilde{J}_\lambda
&=\sum_\mu d^\mu_\lambda(p)\Tilde{J}_\mu.
\end{align}
Then
we see that
$P_i(0)\Tilde{T}(p)P_i(p)|_{V_i(0)}:V_i(0)\to V_i(0)$
and 
$P_i(0)P_i(p)|_{V_i(0)}:V_i(0)\to V_i(0)$ are real-holomorphic.
Equivalently, $\sum_\mu c^\mu_\lambda(p)d^\nu_\mu(p)$ and 
$d^\mu_\lambda(p)$ are real-holomorphic.
By Proposition \ref{proj2}, $P_i(0)P_i(p)|_{V_i(0)}$ or 
the matrix $d^\mu_\lambda(p)$ is invertible,
which implies $c^\mu_\lambda(p)$ is real-holomorphic.
\end{proof}
The matrix $c(p)=(c^\mu_\lambda(p))$ is symmetric.
It is known that if all the matrix elements of the symmetric operator on the finite dimensional vector space are real-holomorphic,
then its eigenvalues and eigenvectors are real-holomorphic. 
(See \cite{Kat})

Hence we have
\begin{prop}
\label{prop:eigen_hol}
The eigenvalues of
$\Tilde{T}(p)$ are real-holomorphic
and coincide with $a_i$ when $p=0$.
The eigenfunctions are also real-holomorphic.
\end{prop}

Summarizing, we obtain the following theorem.
\begin{thm}
\label{thm:eigen_hol}
For each $a_i \in \sigma(\Tilde{H}_0)$, 
there exists $p_i>0$ such that 
for $-p_i<p<p_i$,  
the dimension of the eigenspace whose eigenvalues are
included in $|\zeta-a_i|<W_{\max}(p)$
is equal to
the dimension of the eigenspace of eigenvalue $a_i$.
Moreover
the eigenfunctions and the eigenvalues depend on $p$ 
real-holomorphically.
\end{thm}

If the coupling constants $k_{\alpha}(>0)$ are all rational numbers, 
we can estimate the eigenvalues uniformly. We will explain this below.

Suppose $k_{\alpha} $ are all rational. 
Let $k_{\alpha} =k_{\alpha,\mathrm{num}}/k_{\mathrm{den}}$
such that
$k_{\alpha,\mathrm{num}}$ are integers and $k_{\mathrm{den}}$ is positive integer.
Let $n$ be the minimal positive integer such that $(P|P)\subset\mathbb{Z}/n$.
Then we see that the spectrum of $\Tilde{H}_0$ is uniformly separated.
To be more precise,
if $a,b \in \sigma (\Tilde{H}_0)$ and $a\neq b$,
we have $|a-b|\geq 1/n k_{\mathrm{den}}$.
Hence if we take $p_0$ as
\begin{equation}
\label{p0}
W_{\max}(p_0)=1/4nk_{\mathrm{den}},  
\end{equation}
then
there exists a set of
circles $\Gamma_i$ 
such that for $|p|<p_0$, 
each $\Gamma_i\subset\rho(\Tilde{T}(p))$
contains only one element 
$a_i \in \sigma (\Tilde{H}_0)$ inside it, 
any two circles never cross, 
$\| P_i(p)-P_i(0) \| <1 $,
and every element of $\sigma (\Tilde{T}(p))$ is contained inside of some circle $\Gamma_i $.

Therefore we have
\begin{thm}
Suppose $k_{\alpha} \in \Rat_{>0}$ and 
let $p_0$ be defined in (\ref{p0}).
Then Theorem \ref{thm:eigen_hol} holds for $p_i=p_0$.
All eigenvalues of $\Tilde{T}(p)$ on the $L^2(T,\Delta^2 d\mu)$ space are contained in $\cup _{a \in \sigma (\Tilde{H}_0)} \{ \zeta | \:  |\zeta -a| <W_{\max}(p) \}$ for $-p_0<p<p_0$.
All eigenfunctions are real-holomorphically connected to the eigenfunction of $\Tilde{H}_0$ as $p \rightarrow 0$.
\end{thm}

\section{$A_{N-1}$-cases}
\subsection{}
In section \ref{l2}, we considered the spectrum problem of the gauge-transformed Hamiltonian $\Tilde{T}(p)$ in the $L^2(T,\Delta^2 d\mu)^{W}$ space and show that the perturbation is holomorphic by use of the theory of Kato and Rellich.

On the other hand, it is known that there are commuting family of differential operators (e.g. (\ref{comm}) for the $A_{N-1}$ case) which commute with the Hamiltonian. (\cite{OP,OOS,Has})

In this section, we will investigate the relationship between the functions obtained by applying the projections $P_i(p)$ and the commuting family of differential operators. As a result, we will prove that the perturbation series which is obtained by the algorithmic calculation is not only square-integrable but also holomorphic w.r.t. the variables of the coordinate.

For this purpose, we will consider the spectrum problem in the $C^{\omega }(T)^{W}$-space.

In this section, we consider the $A_{N-1}$ cases.

\subsection{}
We introduce some known result for the $A_{N-1}$ cases. 

We realize the $A_{N-1}$ root system in $\Rea ^N$.
Let $\{ \epsilon_i \} _{i=1, \dots ,N}$ be an orthonormal basis.
The space $\Ha ^{*}$ is defined by $\Ha ^{*} := \{ h=\sum_{i=1}^{N}h_i \epsilon_i | \sum_{i=1}^{N}h_i=0\}$. The simple roots are $\{ \epsilon _i -\epsilon _{i+1} | i=1 ,\dots ,N-1\}$.
We set $x_i=(h|\epsilon _{i})$.

Let us remind the Hamiltonian of the elliptic Calogero-Moser model,
\begin{equation*}
\mathcal{H} :=- \frac{1}{2} \sum_{i=1}^{N} 
\frac{\partial ^{2}}{\partial x_{i}^{2}}
+ 
\beta (\beta  -1)
\sum_{1 \leq i<j \leq N}
\wp ( x_{i}-x_{j}).
\end{equation*}

This system is integrable, \emph{i.e.}, 
there exists sufficiently many commuting operators.
The existence and the explicit expressions are known in (\cite{OP,OOS,Has}) etc.
Here, we exhibit the Hasegawa's expression which will be used in the proof of Proposition \ref{prop:actwelldef}. Later we will discuss the relationship between the expression of Ochiai-Oshima-Sekiguchi (\cite{OOS}) and the one of Hasegawa (\cite{Has}).

Following (\cite{Has}), we set 
\begin{align}
& \hat{\mathcal{H}_i}:=\sum_{|I|=i}\sum_{J\subset I}\frac{
\left( \prod_{j\in J}\beta  \frac{\partial }{\partial x_j} \right) \Theta (x)}{\Theta (x)}
\prod_{j \in I\setminus J}\left(\frac{\partial }{\partial x_j} \right) ,
\label{comm1} \\
& \mathcal{H}_i:=\Theta (x)^{\beta }  \hat{\mathcal{H}_i}\Theta (x)^{-\beta } 
\; \; \; \; (1\leq i\leq N).
\label{comm}
\end{align}
where $\Theta (x):=\prod_{1\leq i<j\leq N}\theta ((x_i -x_j)/2\pi)$, $\theta (x)$ is the theta function defined in section \ref{fns}.

The operators $\hat{\mathcal{H}_i}, \: {\mathcal{H}_i}, \: {\mathcal{H}}$ act on the space of functions which are meromorphic except for the branches along $x_j-x_k \in 2\pi (\Zint + \Zint \tau )$ $(j\neq k)$.

On this space, we have $[ \hat{\mathcal{H}}_{i_1},  \hat{\mathcal{H}}_{i_2}]=[\mathcal{H}_{i_1} , \mathcal{H}_{i_2} ]=0$ $(1 \leq {i_1},{i_2} \leq N)$ and $[\mathcal{H}_i , \mathcal{H}]=0$ $(1 \leq i \leq N)$.

We set $\tilde{\Delta }:= \prod_{1\leq j<k\leq N} \sin ^{\beta}((x_j -x_k)/2)$.
The function $|\tilde{\Delta }|$ is the ground-state of the trigonometric Calogero-Sutherland model (\ref{gstate}), i.e. $|\tilde{\Delta } |=\Delta $.
Although the function $\sin ^{\beta}((x_j -x_k)/2)$ may have branch along $x_j-x_k \in 2\pi \Zint$ $(j\neq k)$, the operators $\tilde{\Delta }^{-1}\mathcal{H} \tilde{\Delta }$ and $\tilde{\Delta }^{-1}\mathcal{H}_i \tilde{\Delta }$ do not have branch points if we rewrite the operators by using the commutation relations for $\tilde{\Delta }$ and $\frac{\partial }{\partial x_i}$.
We define the operators $\tilde{\Delta }^{-1}\mathcal{H} \tilde{\Delta }$ and $\tilde{\Delta }^{-1}\mathcal{H}_i \tilde{\Delta }$ by the form on which the branch points had disappeared and we set $H^{(i)}(p):= \tilde{\Delta }^{-1}\mathcal{H}_i \tilde{\Delta }$ $(1\leq i \leq N)$.
The action of the operator $\tilde{\Delta }^{-1}\mathcal{H} \tilde{\Delta }$ coincide with the one of the operator $\tilde{T}(p)$ for the smooth functions on the real domain except for $x_j=x_k$ $(j\neq k)$.
From this reason, we adopt the notation $T(p)=\tilde{\Delta }^{-1}\mathcal{H} \tilde{\Delta }$. The operators $T(p)$ and $H^{(i)}(p)$ act on the space of meromorphic functions on the complex domain.

The operator $T(p)$ is expressed by some combinations of $H^{(1)}(p)$ and $H^{(2)}(p)$.

\begin{prop} \label{prop:actwelldef}
The operators $T(p)$ and $H^{(i)}(p)$ preserve the space $C^{\omega }(T)^W$, where $T$ is the torus $\mathfrak{h}_{\mathbb{R}}/2\pi{Q}{}^\vee$. ($Q{}^\vee(\simeq Q$): the coroot lattice of type $A_{N-1}$)
\end{prop}
\begin{proof}
Since the operator $T(p)$ is expressed in terms of $H^{(1)}(p)$ and $H^{(2)}(p)$, it is enough to show the $H^{(i)}(p)$ $(1\leq i \leq N)$ cases.
The Weyl group of type $A_{N-1}$ acts on the space of function on $\mathfrak{h}_{\mathbb{R}}$ by the permutation of the variable. We denote the action of $\sigma $ on $f(x)$ by $f(x^{\sigma })$.

Let us recall the definition of $C^{\omega }(T)^W$, i.e.,
\begin{equation}
C^{\omega }(T)^W = \left\{ f(x) \in C^{\omega }({\mathbb{R}}^N) \left|
\begin{array}{l}
 f(x^{\sigma })= f(x) \; (\forall \sigma \in W), \\
 f(x+u\sum_{i=1}^N \epsilon _i)=f(x) \; (\forall u \in {\mathbb{R}}), \\
 f(x+2\pi (\epsilon _i -\epsilon _j))=f(x) \; (1\leq  i,j \leq N)  
\end{array}
 \right. \right\} .
\label{ComegaT}
\end{equation}
Let $f(x)$ be a function in $C^{\omega }(T)^W $. From the definition of the operators $H^{(i)}(p)$ $(1\leq i \leq N)$, the function $\tilde{f}(x):=H^{(i)}(p)f(x)$ satisfies the relations $\tilde{f}(^{\sigma }x)= \tilde{f}(x) \; (\forall \sigma \in W),$ $\tilde{f}(x+u\sum_{j=1}^N \epsilon _j)=\tilde{f}(x) \; (\forall u \in {\mathbb{R}})$, $ \tilde{f}(x+2\pi (\epsilon _{j_1} -\epsilon _{j_2}))=\tilde{f}(x) \; (1\leq  {j_1},{j_2} \leq N)$.

It remains to show the holomorphy of the function $\tilde{f}(x)$ on ${\mathbb{R}}^N$.

The function $(\theta (x) /\sin \pi x)^{\beta}$ is non-zero holomorphic function in ${\mathbb{R}}$, because the function $\theta (x) /\sin \pi x$ is non-zero on ${\mathbb{R}}$ and does not admit any branching points on ${\mathbb{R}}$. From the definition of the operators $\hat{\mathcal{H}}_i$ $(1\leq i \leq N)$,  the function $H^{(i)}(p)f(x)$ does not have poles except for $x_{j_1}-x_{j_2}+2\pi k =0$ $(1\leq {j_1}\neq  {j_2} \leq N, \; k \in \Zint)$ on ${\mathbb{R}}^N$.
If the function $H^{(i)}(p)f(x)$ has a pole along $x_{j_1}-x_{j_2}=0$, the order of the pole is one, but it contradicts to the Weyl group invariance of the function $H^{(i)}(p)f(x)$. Therefore the function $H^{(i)}(p)f(x)$ is holomorphic along $x_{j_1}-x_{j_2}=0$.

The holomorphy along $x_{j_1}-x_{j_2}+2\pi k =0$ $( k \in \Zint \setminus \{ 0 \})$ follows from the periodicity of  $H^{(i)}(p)f(x)$.
\end{proof}
From the commutativity of $\mathcal{H}$ and $\mathcal{H}_i$ we have 
\begin{equation}
[T(p), H^{(i)}(p)]=[H^{(i)}(p), H^{(j)}(p)]=0, \; (1\leq i,j\leq N),
\label{THcomm}
\end{equation}
on the space $C^{\omega }(T)^W $.

By using the formula (\ref{th1}), we can expand the operators $T(p)$ and $H^{(i)}(p)$ $(1\leq i \leq N)$ as the formal power series of operators w.r.t. the parameter $p(=e^{2\pi \sqrt{-1} \tau})$,
\begin{align}
& T(p)= T(0) + \sum_{j=1}^{\infty} T^{\{ j \}}(0) p^j ,\\
& H^{(i)}(p) = H^{(i)}(0) +\sum_{j=1}^{\infty} H^{(i),\{ j \}}(0) p^j .\nonumber
\end{align}
We set $z_i=e^{\sqrt{-1} x_i}$. The operators $T(0) $, $H^{(i)}(0)$, $T^{\{ j \}}(0)$, and $H^{(i),\{ j \}}(0)$ are expressed as the combination of the rational functions of $z_1,z_2, \dots , z_N$ and the polynomials of $\frac{\partial }{\partial z_1}, \dots , \frac{\partial }{\partial z_N}$.

\begin{prop} \label{prop:ffin}
Let $f$ be an element of $\Cplx [P]^W$. Then the functions $T(0) f$, $T^{\{ j \}}(0) f$, $H^{(i)}(0)f$, and $H^{(i),\{ j \}}(0)f$ are elements of  $\Cplx [P]^W$ for $j \in \Zint_{\geq 1}$ and $1\leq i\leq N$.
\end{prop}
\begin{proof}
We prove $H^{(i),\{ j \}}(0) f\in \Cplx [P]^W$. For the other cases, the proofs are similar.

From the definition, the function $H^{(i),\{ j \}}(0) f$ is symmetric and rational with respect to the variables $z_1, \dots ,z_N$.
The possible poles of the rational function $H^{(i),\{ j \}}(0) f$ are $z_k-z_{k'}=0$ $(1\leq k < k' \leq N)$ and the degree of each pole is one, but it contradicts to the Weyl group invariance of the function $H^{(i),\{ j \}}(0) f$.

Therefore the rational function $H^{(i),\{ j \}}(0) f$ does not have any poles, and we have $H^{(i),\{ j \}}(0) \in \Cplx [P]^W$.
\end{proof}

We notice that the operator 
\begin{equation}
T(0)= A_1\Tilde{\mathcal H}_0+A_2
\label{gtHam}
\end{equation}
 is the gauge-transformed Hamiltonian of the trigonometric Calogero-Sutherland model up to constants $A_1, A_2$, where
\begin{equation}
\Tilde{\mathcal H}_0 =\sum_{i=1}^N \left( z_i\frac{\partial}{\partial z_i} \right) ^2+\beta \sum_{i<j}\frac{z_i+z_j}{z_i-z_j}\left( z_i\frac{\partial}{\partial z_i}-  z_j\frac{\partial}{\partial z_j}\right) .
\end{equation}

The joint eigenfunction of the operators $H^{(i)}(0)$ $(1\leq i \leq N)$ is the Jacobi polynomial $J_{\lambda}$ $(\lambda \in P_+)$.
We denote the joint eigenvalue $E^{(i)}_{\lambda }$ by $H^{(i)}(0)J_{\lambda}= E^{(i)}_{\lambda } J_{\lambda}$.

Suppose $\beta >0$. Let $\lambda ,\mu \in P_+$. It is known that the condition $E^{(i)}_{\lambda }=E^{(i)}_{\mu }$ for all $i \in \{ 1,\dots ,N\}$ is equivalent to $\lambda =\mu$. In other words, the joint eigenvalue is non-degenerate.

From now on we will discuss the symmetry (self-adjointness) of the higher commuting Hamiltonians. 
For this purpose, we will discuss the relationship between the expressions of the higher commuting Hamiltonians in (\cite{OOS}) and the ones in (\cite{Has}).

Following (\cite{OOS,OS}), we introduce the operators
\begin{align}
& I_k=\sum_{0\leq j\leq [\frac{k}{2}]}\frac{1}{2^j j! (k-2j)!}  \sum_{\sigma \in W}\:  ^{\sigma }\left( u(x_1-x_2)u(x_3-x_4)  \dots  \right. \label{OOScommH} \\
& \left. \dots u(x_{2j-1}-x_{2j}) \frac{\partial}{\partial x_{2j+1}}  \frac{\partial}{\partial x_{2j+2}}\dots \frac{\partial}{\partial x_{k}} \right), \nonumber 
\end{align}
where $k=1,\dots ,N$, $W$ is the Weyl group of $A_{N-1}$-type ($N$-th symmetric group),  $^{\sigma }(f(x_1, \dots ,x_N))= f(x_{\sigma ^{-1} (1)}, \dots ,x_{\sigma ^{-1} (N)})$ for $\sigma \in W$, and $u(x)=\beta (\beta -1) \wp(x)$.
The domains of the operators $I_k$ $(k=1,\dots ,N)$ are the same as the ones of $\mathcal{H}_k$ $(k=1,\dots ,N)$.

By a straightforward calculation, we have $\mathcal{H}_3=I_3 +C I_1$ for some constant $C$.
Applying  Theorem 5.2. in (\cite{OS}), we obtain that the operators $\mathcal{H}_k$ $(k=1, \dots ,N)$ are expressed as the polynomial of $I_1, I_2, \dots , I_N$.

Let $\Rea [I_1, I_2, \dots , I_N]$ be a polynomial ring generated by $I_1, I_2, \dots , I_N$ and $\varsigma $ be an involution on $\Rea [I_1, I_2, \dots , I_N]$ such that $^{\varsigma }F(x_1, \dots ,x_N)=F(-x_1,\dots , -x_N)$. Then $^{\varsigma }I_k=(-1)^k I_k$ and $^{\varsigma } \mathcal{H}_k= (-1)^k \mathcal{H}_k$.
Hence $\mathcal{H}_k$ admit the expansion,
\begin{equation}
\mathcal{H}_k =\sum_{j_1\leq \dots \leq  j_m}c_{j_1, \dots ,j_m}I_{j_1}  \dots I_{j_m},
\label{Hkeo}
\end{equation}
where $c_{j_1, \dots ,j_m} \in \Rea$ and if $k-(j_1+ \dots +  j_m) \not \in 2\Zint_{\geq 0}$ then $c_{j_1, \dots ,j_m} =0$.

From the similar discussion, the operators $I_k$ 
admit the expansion,
\begin{equation}
I_k =\sum_{j_1\leq \dots \leq  j_m}\tilde{c}_{j_1, \dots ,j_m}\mathcal{H}_{j_1}  \dots \mathcal{H}_{j_m},
\label{Ikeo}
\end{equation}
where $\tilde{c}_{j_1, \dots ,j_m} \in \Rea$ and if $k-(j_1+ \dots +  j_m) \not \in 2\Zint_{\geq 0}$ then $\tilde{c}_{j_1, \dots ,j_m} =0$.

\begin{lemma} \label{symlem}
We suppose $\beta > N$.
For $f,g \in C^{\omega }(T)^W $, we have $(H^{(k)}(p)f,g)_{\Delta}=(-1)^k (f,H^{(k)}(p)g)_{\Delta}$ $(1\leq k \leq N)$.
\end{lemma}

\begin{proof}
Since $(f,g)_{\Delta}= \int_T \overline{f(x)} g(x) |\tilde{\Delta}|^2 d\mu $ and $H^{(k)}(p)= \tilde{\Delta}^{-1} \mathcal{H}_k \tilde{\Delta}$, it is enough to show $\int_T   \left(\overline{\mathcal{H}_k (f(x)|\tilde{\Delta}|)} \right) g(x)|\tilde{\Delta}|  d\mu = (-1)^k \int_T \overline{ f(x)|\tilde{\Delta}|}   \left( \mathcal{H}_k (g(x)|\tilde{\Delta}|)\right) d\mu $.
We have the equality $\int_{T}h(x) d\mu = A \int_{0\leq x_1, \dots ,x_N \leq 2\pi N}h(x) dx_1 \dots dx_N$ for some non-zero constant $A$, which follows from the correspondence between the integration of the $sl_n$-invariant function and the one of the $gl_n$. From this equality,  the property (\ref{Hkeo}), and the commutativity $[I_{j_1}, I_{j_2}]=0$ $(1\leq j_1 ,j_2 \leq N)$, if we show
\begin{equation}
\int_D  \left(  \overline{I_k( f(x)|\tilde{\Delta}|)} \right) g(x)|\tilde{\Delta}| dx = (-1)^k \int_D  \overline{f(x)|\tilde{\Delta}|}  \left( I_k (g(x)|\tilde{\Delta}|)\right)dx,
\label{Isymm}
\end{equation}
where $D=\{(x_1, \dots ,x_N) \in \Rea^N | 0\leq x_1, \dots ,x_N \leq 2\pi N \}$ and $dx=dx_1dx_2\dots dx_N$,
then we obtain Lemma \ref{symlem}.

If $\beta > N$ then the functions $f(x) |\tilde{\Delta}|, g(x)|\tilde{\Delta}|$ are $C^N$-class.

From the definition of $C^{\omega }(T)^W $ (\ref{ComegaT}), we have the periodicity $f(x_1, \dots,  x_l+2\pi N, \dots, x_N)=f(x_1, \dots,  x_l, \dots, x_N)$ $(1\leq l \leq N)$ for $f(x_1, \dots , x_N) \in C^{\omega }(T)^W $.

We set $\tilde{f}(x) =f(x)|\tilde{\Delta}|$ and $\tilde{g}(x) =g(x)|\tilde{\Delta}|$.
The functions $\tilde{f}(x) $, $\tilde{g}(x)$ are smooth on $\Rea ^N$ except for $x_i-x_j+2\pi k =0$ $(1\leq i\neq j \leq N, \; k\in \Zint )$.
The behaviors of the functions $\tilde{f}(x) $, $\tilde{g}(x)$ around $x_i-x_j+2\pi k =0$ are $O(|x_i-x_j|^{\beta})$, i.e. $\frac{\tilde{f}(x)}{|x_i-x_j|^{\beta}}$ and $\frac{\tilde{g}(x)}{|x_i-x_j|^{\beta}}$ are bounded around $x_i-x_j+2\pi k =0$.

From the expression of $I_k$ (\ref{OOScommH}), if we show 
\begin{align}
& \int _{D} \left( u(x_1-x_2)u(x_3-x_4)  \dots   u(x_{2j-1}-x_{2j}) \frac{\partial}{\partial x_{2j+1}} \dots \frac{\partial}{\partial x_{k}} \tilde{f}(x) \right) \tilde{g}(x) dx  \label{partint} \\
& = \int _{D} \tilde{f}(x) \left( u(x_1-x_2)u(x_3-x_4)  \dots   u(x_{2j-1}-x_{2j}) \frac{\partial}{\partial x_{2j+1}} \dots \frac{\partial}{\partial x_{k}} \tilde{g}(x) \right) dx, \nonumber
\end{align}
for all $j$ s.t. $0\leq j \leq [\frac{k}{2}]$, then
we obtain (\ref{Isymm}) and Lemma \ref{symlem}.

The number $\beta$ satisfies $\beta >N\geq 2$.
Though the function $u(x_{2l-1} - x_{2l})= \beta (\beta -1) \wp ( x_{2l-1} - x_{2l})$ $(l=1,\dots ,j) $ have double pole along $x_{2l-1} - x_{2l}+2\pi k=0$ $(k\in \Zint )$, the integrands of (\ref{partint}) are bounded around $x_{2l-1} - x_{2l}+2\pi k=0$ from the properties $\tilde{f}(x) = O(|x_{2l-1} - x_{2l}|^{2})$ and $\tilde{g}(x) = O(|x_{2l-1} - x_{2l}|^{2})$ around $x_{2l-1} - x_{2l}+2\pi k=0$.
Hence the singularities along $x_{2l-1} - x_{2l}+2\pi k=0$ $(l=1,\dots ,j, \; k\in \Zint )$ do not affect the integration.

Since the integrands of (\ref{partint}) are continuous,
we can replace the range of integration of both sides of (\ref{partint}) with $D'$, where $D'= \{ (x_1,\dots x_N) \in D |  x_{2l-1}- x_{2l} +2\pi k \neq 0  \; (l=1,\dots ,j, \;  k\in \Zint ) \}$.

It is obvious that 
\begin{align}
& \int _{D'} \left( u(x_1-x_2)u(x_3-x_4)  \dots   u(x_{2j-1}-x_{2j}) \frac{\partial}{\partial x_{2j+1}} \dots \frac{\partial}{\partial x_{k}} \tilde{f}(x) \right) \tilde{g}(x) dx  \label{intDprime}  \\
& = \int _{D'} \left( \frac{\partial}{\partial x_{2j+1}} \dots \frac{\partial}{\partial x_{k}} u(x_1-x_2)u(x_3-x_4)  \dots   u(x_{2j-1}-x_{2j}) \tilde{f}(x) \right) \tilde{g}(x) dx. \nonumber
\end{align}
By applying the integration by parts repeatedly, we find that the r.h.s. of (\ref{intDprime}) is equal to 
\begin{align}
& (-1)^{k-2j} \int _{D'} u(x_1-x_2)u(x_3-x_4)  \dots   u(x_{2j-1}-x_{2j}) \tilde{f}(x) \frac{\partial}{\partial x_{2j+1}} \dots \frac{\partial}{\partial x_{k}}\tilde{g}(x) dx.\nonumber
\end{align}
Here we used the periodicities on $x_l \rightarrow x_l+2\pi N$ ($l=2j+1, \dots ,k$).

Hence we obtain (\ref{partint}) and Lemma \ref{symlem}.

\end{proof}

\begin{prop} \label{symprop}
We suppose $\beta \geq 0$.
For $f,g \in C^{\omega }(T)^W $,  we have $(H^{(k)}(p)f,g)_{\Delta}=(-1)^k (f,H^{(k)}(p)g)_{\Delta}$ $(1\leq k \leq N)$. In other words, the operators $(\sqrt{-1})^k H^{(k)}(p)$ are symmetric on the space $C^{\omega }(T)^W $.
\end{prop}
\begin{proof}
It is trivial for the $\beta =0$ case.

We assume $\beta >0$.

Let $f (x) \in C^{\omega }(T)^W $. Then $H^{(k)}(p)f(x) $ is a polynomial in the parameter $\beta $ of degree at most $k$ and $H^{(k)}(p)f(x) \in C^{\omega }(T)^W $.

We set $H^{(k)}(p)f(x)  = \sum_{j=0}^k f_j(x)\beta ^j$. Then $f _j(x) \in C^{\omega }(T)^W $ $(0\leq j\leq k)$, because $H^{(k)}(p)f(x) \in C^{\omega }(T)^W $ for all $\beta$.
For $f(x)$ and $H^{(k)}(p)f(x)$, we set $\widetilde {f(x)}= \overline{f(x)}$ and $\widetilde {H^{(k)}(p) f(x)}= \sum_{j=0}^k \overline{f_j(x)}\beta ^j$.

We fix the functions $f(x) ,g(x) \in C^{\omega }(T)^W $.
It is enough to show that the equations $(H^{(k)}(p)f,g)_{\Delta}-(-1)^k (f,H^{(k)}(p)g)_{\Delta}=0$ hold for $\beta >0$ and $1\leq k \leq N$.

We set
\begin{align*}
& T'= \{ (x_1,\dots ,x_N)\in \Rea ^N| \sum_{i=1}^N x_i=0, \; 0\leq x_i-x_j\leq 2\pi \: (1\leq i<j\leq N)\}, \\ 
& \stackrel{\circ}{T'}= \{ (x_1,\dots ,x_N)\in \Rea ^N|  \sum_{i=1}^N x_i=0, \; 0< x_i-x_j< 2\pi \: (1\leq i<j\leq N)\},\\
& \stackrel{\ast}{h} \! (x)=  \left(\overline{H^{(k)}(p)f(x)} \right) g(x)|\tilde{\Delta}|^2 - (-1)^k \overline{f(x)} \left( H^{(k)}(p) g(x)\right) |\tilde{\Delta}| ^2,\\
& h(x)=  \left(\widetilde{H^{(k)}(p)f(x)}\right) g(x)\tilde{\Delta}^2 - (-1)^k \widetilde{f(x)} \left( H^{(k)}(p) g(x)\right)\tilde{\Delta}^2 .
\end{align*}

Then the equation $(H^{(k)}(p)f,g)_{\Delta}-(-1)^k (f,H^{(k)}(p)g)_{\Delta}=0$ is equivalent to $\int_T \!  \stackrel{\ast}{h} \! (x)d\mu =0$, where $T=\mathfrak{h}_{\mathbb{R}}/2\pi{Q}{}^\vee$. From the equation $ \frac{1}{N!} \int_T \! \stackrel{\ast}{h}\! (x)d\mu = \int_{T'} \! \stackrel{\ast}{h}\! (x)d\mu $, it is sufficient to show $\int_{T'}\!  \stackrel{\ast}{h} \! (x)d\mu =0$.

We have $h(x)=\stackrel{\ast}{h}\! (x)$ on the domain $\stackrel{\circ}{T'}$, because $\sin((x_i -x_j)/2)>0$ on $\stackrel{\circ}{T'}$ for $i<j$ and the branch of the function $\sin^{\beta}((x_i -x_j)/2)$ is chosen to be a positive real number.
For $\beta \in \Cplx$, the branch of the function $\sin^{\beta}((x_i -x_j)/2)$ $(i<j)$ is canonically chosen by the relation $a^{\beta } =\exp (\beta \log a)$ for $a= \sin ((x_i -x_j)/2)>0$.
Hence it is sufficient to show the equation
$\int_{T'} h(x)d\mu =0$ for $\beta >0$.

From Lemma \ref{symlem}, $\int_{T'} h(x)d\mu =\int_{T'} \stackrel{\ast}{h}\! (x)d\mu =0$ holds for $\beta >N$.

From Proposition \ref{prop:actwelldef}, we have $ H^{(k)}(p) f(x) \in C^{\omega }(T)^W $ when $f(x) \in C^{\omega }(T)^W $ and $\beta \in \Cplx$. Hence the integral $\int_{T'} h(x)d\mu $ is well-defined if $\mbox{Re}\beta >0$.

We fix $\beta_0 $ ($\mbox{Re}\beta _0 >0$). Since the function $h(x)$ is holomorphic in $\beta$ and the functions $h(x)$ and $\frac{\partial}{\partial \beta }h(x)$ are uniformly bounded in $(x,\beta ) \in T' \times \{ \beta '| \: |\beta '-\beta_0|<\epsilon \}$ for some $\epsilon \in \Rea _{>0}$, the integral  $\int_{T'} h(x)d\mu $ is also holomorphic at $\beta=\beta_0 $ $(\mbox{Re}\beta _0 >0)$ by the Lebesgue's theorem.

By the identity theorem, the equation $\int_{T'} h(x)d\mu =0$ holds for $\beta$ s.t. $\mbox{Re}\beta >0$.
Therefore we obtain the proposition.
\end{proof}

\subsection{Perturbation} \label{secpert}

We start with the general proposition related to the perturbation method.

\begin{prop} \label{pert}
Let $\{ v_1, v_2, \dots \}$ be linearly independent vectors in a vector space $V$ over $\Rea$.
Let $H_i^{\{ k \}}$ $(k\in \Zint _{\geq 0},\;  i=1, \dots ,N)$ be linear operators on $V$ such that $
H_i^{\{ k \}} v_j=\sum_{j': \mbox{\scriptsize{finite}}} d^{\{ k \},i}_{j,j'}v_{j'} \; \; $ for all $i,j,k$. We assume that there exists 
$ E^{\{ 0\},i }_{j}\in \Rea$ such that
$H_{i}^{\{ 0\} } v_j= E^{\{ 0\} ,i}_{j}v_j$ for all $i,j$ and if $j_1 \neq j_2$ then there exists some $i$ such that $ E^{\{ 0\} ,i}_{j_1} \neq E^{\{ 0\} ,i}_{j_2}$.
Let $( \: , \: )$ be an inner product on $V$ such that $(v_i,v_j) =\delta _{i,j}$.
Let $H_i (p):= \sum_{k=0}^{\infty} H_i ^{\{ k \}} p^k$ be formal power series of the linear operators and assume $[H_{i_1}(p), H_{i_2}(p)] =0$ for all $i_1,i_2 \in \{1,\dots ,N\}$ as the formal power series of $p$.
Then there exists formal power series of vectors 
\begin{equation}
v_j(p)= v_j+\sum_{k=1}^{\infty} \sum_{j':  \mbox{\scriptsize{finite}}} c^{\{ k \} }_{j,j'}v_{j'}p^k,
\label{fpseries}
\end{equation}
such that $H_i(p)v_j(p)= E_{j}^{i}(p)v_j(p)$ and $(v_j(p), v_j(p))=1$, where $E_{j}^{i}(p)= E_{j}^{\{ 0\} ,i} +\sum_{k=1}^{\infty}E_{j}^{\{ k\} ,i} p^k$ is a formal power series on $p$ and the equalities hold as the formal power series of $p$.

For each $j$, the normalized formal power series of the joint eigenfunction of the form (\ref{fpseries}) is unique. 
\end{prop}
\begin{proof}
We introduce variables $w_1, \dots ,w_N$ and set
\[
H(w,p):=\sum_{i=1}^N w_iH_i(p), \; \;
\sum_{j'}d^{\{ k \}}_{j,j'}(w)v_{j'}:=\sum_{i=1}^N w_iH_i^{\{ k \}} v_j,
\]
\[
v_j(p):= v_j+\sum_{k=1}^{\infty } \sum_{j'} c^{\{ k \}}_{j,j'}v_{j'}p^k,
\]
\[ 
E_j (w,p):=\sum_{k=0}^{\infty } E_j^{\{ k \}}(w)p^k=\sum_{i=1}^N  E_{j}^{\{ 0\} ,i}w_i+ \sum_{k=1}^{\infty } \sum_{i=1}^N  E_{j}^{\{ k\} ,i}w_ip^k.
\]
The numbers $d^{\{ k \}}_{j,j'}(w)$ and $E_{j}^{\{ 0\} ,i}$ are given in advance.
We will investigate the conditions for the coefficients of the formal power series $v_j(p)$ and $E_j{(w,p)}$ satisfying the following relations
\begin{equation}
 H(w,p)v_j(p)=E_j {(w,p)}v_j(p), \; \; 
 (v_j(p), v_j(p))=1. 
\label{pereq}
\end{equation}
We set $c_{j,j'}^{\{ 0\} }=\delta_{j,j'}$
By comparing the coefficients of $v_{j'}p^k$, we obtain that the conditions (\ref{pereq}) are equivalent to the following relations,
\begin{align}
& c_{j,j'}^{\{ k \}}= \frac{\sum_{k'=1}^{k} 
( \sum _{j''}c_{j,j''}^{\{ k-k' \} }d_{j'',j'}^{\{ k' \}}(w))-\sum_{k'=1}^{k-1} c_{j,j'}^{\{ k-k' \} }E_{j}^{\{ k' \}}(w) }
{E_j^{\{ 0 \}}(w)-E_{j'}^{\{ 0 \}}(w)} ,\; \;(j' \neq j) \label{nondeg} \\
& c_{j,j}^{\{ k \}}=-\frac{1}{2} \left(\sum_{k'=1}^{k-1} \sum_{j'} c_{j,j'}^{\{ k' \}}c_{j,j'}^{\{ k-k' \} } \right) ,\\
& E_j^{\{ k \}} (w)= \sum_{k'=1}^{k}\sum_{j'}c_{j,j'}^{\{ k-k' \} }d_{j',j}^{\{ k' \}}(w)- \sum_{k'=1}^{k-1} c_{j,j}^{\{ k-k' \} } E_j^{\{ k' \}}(w).\label{nondeg2}
\end{align}
We remark that the denominator of (\ref{nondeg}) is non--zero by the non--degeneracy condition.
The numbers  $c_{j,j'}^{\{ k \}}$ and $E_j^{\{ k \}}(w)$ are determined recursively and they exist uniquely.
We have recursively that for each $j$ and $k$, $\# \{ j' | \: c_{j,j'}^{\{ k \}} \neq 0 \}$ is finite and the summations in (\ref{nondeg} - \ref{nondeg2}) on the parameters $j'$ and $j''$ are indeed the finite summations.

At this stage, the apparent expression of the coefficients $c_{j,j'}^{\{ k \}}$ may depend on $w$. 
We will show that the coefficients $c_{j,j'}^{\{ k \}}$ do not depend on $w$. We denote $v_j(p)$ by $v_j(w,p)$.

From the commutativity of $H(w,p)$ and $H(w',p)$ we have
\[
H(w,p) (H(w',p) v_j(w,p))=E_j{(w,p)}(H(w',p) v_j(w,p)).
\]
Since the vector $H(w',p) v_j(w,p)$ admits the expansion $H(w',p) v_j(w,p) =E^{\{ 0\} }_j(w') v_j + O(p)$, we obtain the following relation from the uniqueness of the formal eigenvector.
\[
\frac{H(w',p)v_j(w,p)}{\sqrt{f(w,w',p)}} = v_j(w,p).
\]
where $ f(w,w',p):=(H(w',p)v_j(w,p), H(w',p)v_j(w,p))$ and $1/\sqrt{f(w,w',p)}$ is regarded as a formal power series on $p$ from the formula $(a_0^2+\sum_{k=1}^{\infty} a_kp^k)^{-1/2}=a_0^{-1}\left( \sum_{n=0}^{\infty} \left(\!  \begin{array}{c} -1/2 \\ n \end{array} \! \right) \left( a_0^{-2} \sum_{k=1}^{\infty} a_k p^k \right) ^n \right)$.
Therefore we have $H(w',p)v_j(w,p)= \sqrt{f(w,w',p)} v_j(w,p)$. On the other hand we have $H(w',p)v_j(w',p) = E_j(w',p)v_j(w',p)$.
By the uniqueness of the formal eigenvector whose leading term is $v_j$, we have $v_j(w,p)=v_j(w',p)$. Therefore the coefficients $c_{j,j'}^{\{ k \}}$ do not depend on $w$.

From (\ref{nondeg2}), we obtain recursively that the coefficients of the formal eigenvalue $E_j^{\{ k \} }(w)$ are linear in $w_1, \dots ,w_N$.
Therefore the numbers $E_{j}^{\{ k \} ,i}$ $(k\in \Zint_{\geq 1})$ are determined appropriately.
\end{proof}

\begin{prop}
Proposition \ref{pert} is applicable for the $A_{N-1}$-type elliptic Calogero Moser model by the following correspondence,
\begin{eqnarray*}
H_i(p) & \Leftrightarrow & \mbox{The commuting differential operator } H^{(i)}(p), \\
v_j & \Leftrightarrow & \mbox{The normalized Jacobi polynomial } \Tilde{J}_{\lambda} .
\end{eqnarray*}
\end{prop}
\begin{proof}
The finiteness of the summation  $H_i^{\{ k \}} v_j=\sum_{j'} d^{\{ k \},i}_{j,j'}v_{j'} \; \; $ follows from Proposition \ref{prop:ffin} and the fact that the Jack polynomial forms a basis of $\Cplx [P]^W$. 

The non-degeneracy of the joint eigenvalues $E^{\{0\} ,i}_j$ follows from the non-degeneracy of the joint eigenvalue of the Jack polynomial.
\end{proof}

Summarizing, we have the algorithm of computing the ``formal'' eigenvalues and ``formal'' eigenfunctions of the elliptic Calogero-Moser model of $A_{N-1}$-type by using the Jacobi polynomial.
In the next subsection, we will discuss the convergence.

\subsection{Analyticity and the higher commuting operators}
In this subsection, we will consider the spectral problem in the $C^{\omega }(T)^{W}$-space for the $A_{N-1}$ elliptic Calogero-Moser model. We assume $\beta >1$. Since $T$ is compact, we have $C^{\omega }(T)^{W} \subset L^2(T,\Delta^2d\mu)^{W}$.

We will show the holomorphy of the eigenfunctions which we have found on the $  L^2(T,\Delta^2d\mu)^{W}$ space in section \ref{sec:pertL2}.
After having the holomorphy of the eigenfunctions, we will justify the convergence and the holomorphy of the joint eigenfunctions of the higher commuting operators obtained by the algorithmic calculation, which we have explained in section \ref{secpert}.

For this purpose, we need the following propositions.
\begin{prop} \label{analthm} 
For each eigenvalue $a_i\in  \sigma (\Tilde{H}_0)$ and eigenfunction $\Tilde{J}_{\lambda}$ of the Hamiltonian $\Tilde{H}_0$ of the trigonometric model such that $\Tilde{H}_0\Tilde{J}_{\lambda}=a_i\Tilde{J}_{\lambda}$, 
there exists a positive number $p_i$ such that 
the function $P_i(p) \Tilde{J}_{\lambda} $ is holomorphic in $(x_1, \dots x_N,p)$
on the set $B_{p_i}$, where the operator $P_i(p)$ is a projection on the Hilbert space $L^2(T,\Delta^2d\mu)^{W}$ which was defined in section \ref{l2} and 
\begin{equation}
B_{\epsilon }=\{ (x_1,\dots ,x_N, p) \in \Cplx ^N \times \Rea | \: |\mbox{Im} x_j|<\epsilon \: (j=1,\dots ,N), \; -\epsilon<p<\epsilon \}.
\end{equation}
\end{prop}

\begin{prop} \label{prop:commHP}
For all eigenvalue $a_i\in  \sigma (\Tilde{H}_0)$ and Jacobi polynomial $\Tilde{J}_{\lambda} $, we have
\begin{equation}
H^{(j)}(p)P_i(p)\Tilde{J}_{\lambda}=P_i(p)H^{(j)}(p)\Tilde{J}_{\lambda}, \; \; (j=1,\dots ,N), 
\end{equation}
when $|p|$ is sufficiently small.
\end{prop}

We will prove Propositions \ref{analthm} and \ref{prop:commHP} in the next section.

\begin{rmk}
For the $A_1$ and $\beta \in \Zint_{>1}$ cases, and the $A_2$ and $\beta =2$ case, Proposition \ref{analthm} is obvious from the construction of the eigenfunctions via the Bethe Ansatz method (\cite{Tak}).
\end{rmk}

We fix the eigenvalue $a_i \in \sigma (\Tilde{H}_0)$.
From Propositions \ref{prop:multH0}, \ref{analthm}, and \ref{prop:commHP},
if $|p|$ is sufficiently small then
the operators $H^{(j)}(p)$ act on the finite dimensional space $V_i(p)$, where
\begin{equation}
V_i(p)=\sum_{\lambda| \Tilde{H}_0\Tilde{J}_{\lambda}=a_i\Tilde{J}_{\lambda}} \Cplx P_i(p) \Tilde{J}_{\lambda}.
\end{equation}
and we have $V_i(p) \subset C^{\omega}(T)^{W}$.

From Proposition \ref{symprop}, the higher commuting operators $(\sqrt{-1})^j H^{(j)}(p)$ $(j=1,\dots ,N)$ are symmetric both on the space $C^{\omega}(T)^{W}$ and the finite dimensional space $V_i(p)$.
The joint eigenvalues are real-holomorphic w.r.t the parameter $p$ and the operators $(\sqrt{-1})^j H^{(j)}(p)$ are simultaneously diagonalizable in the space $V_i(p)$ if $|p|$ is sufficiently small and $p\in \Rea$.
The joint eigenfunctions are holomorphic on the domain $B_{\epsilon}$ for sufficiently small $\epsilon \in \Rea _{>0}$.

Therefore the joint eigenfunction of $H^{(1)}(p), \dots ,H^{(N)}(p)$ admits the holomorphic expansion in the variable $p$.

Since the joint eigenvalues of $H^{(1)}(0), \dots ,H^{(N)}(0)$ are distinct, the expansion is unique up to the normalization. (See section \ref{secpert}.)
Hence the perturbation series which is obtained by the method introduced in section \ref{secpert} converges holomorphically and 
coincides with the eigenvalue and eigenfunction which is obtained by diagonalizing the finite dimensional space $V_i(p)$.
Summarizing, we have
\begin{thm}
For the $A_{N-1}$ and $\beta >1$ cases, the perturbation expansion of the commuting operators $H^{(1)}(p), \dots ,H^{(N)}(p)$ which is performed in section \ref{secpert} converges holomorphically and defines the eigenfunction which is holomorphic when $|Im x_j |$ $(j=1,\dots ,N)$ and $|p|$ $(p \in \Rea )$ are sufficiently small.
The joint eigenvalue is holomorphic in the parameter $p(\ll 1)$.
\end{thm}

\begin{rmk}
It was pointed out by Prof. T. Oshima that the real-holomorphy of the square-integrable eigenfunction $\psi(x)$ (i.e. $ T(p)\psi(x)=E(p)\psi(x), \; \psi(x) \in L^2(T,\Delta ^2d\mu)^{W}$) is also obtained by the following argument.

From the ellipticity of the operator $T(p)$ and the Weyl's lemma, we have the real-holomorphy of the eigenfunction $\psi(x)$ on the domain $\dot{T}=T\setminus T'$, where $T':=\{ (x_1,\dots ,x_N) \in T | \exists (i\neq j), x_i=x_j \}$.
Next we consider the analytic continuation of the function $\psi(x)$.
The equation $T(p)\psi(x)=E(p)\psi(x)$ has regular singularities along $x_i=x_j \: ( i\neq j)$, and the exponents at the singularity are $(0, -2\beta -1)$. It follows that the function $\psi(x)$ is holomorphic along $x_i-x_j=0$ from the property $\psi(x)\in L^2$.
Hence we have the real-holomorphy of $\psi(x)$ on $T$.
\end{rmk}

\section{Proof of Propositions \ref{analthm} and \ref{prop:commHP}}

In this section, we assume that the root system is of the $A_{N-1}$-type.

For $\lambda \in P_+$ and $j=1,\dots, N-1$, we set $m_{\lambda }= \sum_{\mu \in W\lambda}e^{\langle \mu, h\rangle}$ and $e_{\Lambda _{j}}=m_{\Lambda _{j}}$, where $\Lambda _{j}$ is the $j$-th fundamental weight.
For $\lambda= \sum_{j=1}^{l} \Lambda_{i_j}$ $( l \in \Zint_{\geq 0}, \; i_j \in \{ 1,\dots ,N-1\}$ $(j=1,\dots ,l))$, we set $\tilde{e}_{\lambda }=\prod_{j=1}^l e_{\Lambda_{i_j}}$.
Then we have $\tilde{e}_{\lambda }= e_{\lambda '}$ on $\Ha ^{*}$, where $e_{\mu }$ is the elementary symmetric function for the partition $\mu$ defined in the Macdonald's book (\cite{Mac}, p.20) and $\lambda '$ is the conjugate of $\lambda $.

We set 
$$
t(x,p):=\sum_{k=1}^{\infty}t_k(x)p^k:=\wp (x)
-\frac{1}{4\sin^2 (x/2)}+\frac{1}{12},
$$
which converges uniformly on a strip around $\mathbb{R}\times [-\epsilon ,\epsilon ]$ for $0\leq \epsilon <1$.

From the formula (\ref{wpth}), we have
\begin{equation}
t_k(x)=-2\sum _{j|k} j(\cos jx -1).
\label{tkx}
\end{equation}
Here, $j|k$ means that the positive integer $j$ is a divisor of $k$.

\begin{lemma} \label{wplem}
For a real number $c$ such that $c>1$, there exists a positive number $a'$ such that $|t_k(x)|<a'c^k$ for all $x\in \Rea$.
\end{lemma}
\begin{proof}
Let $t_k$ be the sum of all divisors of $k$.
By the formula (\ref{tkx}), we have $|t_k(x)|\leq 4t_k<4k^2$.
Since the convergence radius of the series $\sum k^2 p^k$ is equal to $1$, the convergence radius of the series $\sum t_k p^k$ is equal to or less than $1$. Therefore we have the lemma.
\end{proof}
The operator $W(p)$ defined in (\ref{eq:Wp}) has an expansion in terms of $p$ given by
$W(p)=\sum_{k=1}^{\infty} T^{(k)}p^{k}$, where
$T^{(k)}$
is the operator of multiplication by the function
$T^{(k)}(h):=\sum_{\alpha \in R_+}
\beta(\beta-1)
t_k(\langle \alpha, h \rangle )$.
For each $p \in (1,-1)$, the series $\sum_{k=1}^{\infty} T^{(k)}p^{k}$ converges uniformly on $\mathfrak{h}_{\mathbb{R}}$:
\begin{prop} \label{tkineq}
For a real number $c$ such that $c>1$, there exists a positive number $a$ such that
 $\| T^{(k)} \| \leq ac^k$.
\end{prop}
\begin{proof}
It follows from Lemma \ref{wplem} and the inequality 
\begin{equation}
\int_{T} |T^{(k)} f|^2 \Delta ^2 d\mu \leq \sup_{h \in T} |T^{(k)}(h)|^2 \int_{T} |f|^2 \Delta ^2 d\mu
\end{equation}
\end{proof}

\begin{lemma} \label{anallem1}
The function $\sum_{\alpha \in \Delta_{+}}t_k(\langle \alpha, h\rangle )$ admits the expansion,
\begin{equation}
\sum_{\alpha \in \Delta_{+}}t_k(\langle \alpha, h\rangle )=\sum_{\mu \in Q \cap P_+ ,|\mu |\leq \sqrt{2}k}c_{\mu}m_{\mu}.
\end{equation}
\end{lemma}
\begin{proof}
From the formula (\ref{tkx}), we have 
$$ \sum_{\alpha \in \Delta_{+}}t_k(\langle \alpha, h\rangle ) = -\sum_{\alpha \in \Delta_{+}} (\sum_{j|k}j(e^{j\langle \alpha, h\rangle }+e^{-j\langle \alpha, h\rangle }-2)
= -\sum_{j|k} 2j(m_{j\theta }-1),
$$
where $\theta$ is the highest root of the root system $A_{N-1}$.
Since $|\theta |=\sqrt{2}$ and $j\theta \in Q \cap P_+$, we have the lemma.
\end{proof}

\begin{suble} \label{sublemma1}
Let $J_{\lambda}$ be the $A_{N-1}$-Jacobi polynomial. We have
\begin{equation}
J_{\lambda }e_{\Lambda_r}=\sum_{\nu \in P_+, \nu-\lambda \in \{ w\Lambda_r | w \in W\} }\bar{c}_{\nu} J_{\nu},
\end{equation}
for some constants $\bar{c}_{\nu}$.
\end{suble}
\begin{proof}
This follows from the Pieri formula (\cite{Mac}, p.332 and section VI.10.).
\end{proof}

\begin{suble} \label{sublemma2}
Let $l$ be a positive integer.
Assume $i_j \in \{ 1,\dots ,N-1\}$ and $w_j \in W$, $(j=1,\dots ,l)$. We have $|\sum_{j=1}^{l} w_j(\Lambda _{i_j})| \leq |\sum_{j=1}^{l} \Lambda _{i_j}|$.
\end{suble}
\begin{proof}
It is sufficient to show $(\lambda +\mu , \lambda +\mu )\geq (\lambda +w(\mu) , \lambda +w(\mu) )$ for $\lambda ,\mu \in P_+$ and $w\in W$. This inequality is equivalent to $(\lambda , \mu-w(\mu))\geq 0$. From the property $ \mu-w(\mu) \in Q_+$, we have $(\lambda , \mu-w(\mu))\geq 0$.
\end{proof}

\begin{suble} \label{sublemma4}
If $\lambda, \mu \in P_+$ and $\lambda - \mu \in Q_+$, then we have $|\lambda | \geq |\mu|$.
\end{suble}
\begin{proof}
Immediate from the equality $(\lambda ,\lambda ) -(\mu ,\mu)=(\lambda -\mu,\lambda+\mu)$. 
\end{proof}

\begin{suble}$($\cite{Mac}, p.20$)$ \label{sublemma3}
The monomial symmetric function $m_{\lambda}$ has the expansion 
\begin{equation}
m_{\lambda} = \tilde{e}_{\lambda}+ \sum_{\nu \in P_+, \lambda -\nu \in Q_+\setminus \{0\}}\check{c}_{\nu} \tilde{e}_{\nu},
\end{equation}
for some constants $\check{c}_{\nu}$.
\end{suble}

\begin{lemma} \label{anallem2}
We have the expansion,
\begin{equation}
J_{\lambda }m_{\mu }=\sum_{\nu \in P_+, |\nu-\lambda |\leq |\mu| }\bar{c}_{\nu} J_{\nu},
\end{equation}
for some constants $\bar{c}_{\nu} $.
\end{lemma}
\begin{proof}
First, we expand $m_{\mu }$ by using Sublemma \ref{sublemma3}. Then $J_{\lambda }m_{\mu }$ is expressed as the linear combination of $J_{\lambda }\tilde{e}_{\nu }$, where $\nu \in P_+$ and $\mu - \nu \in Q_+$. 
We set $\nu =\sum_{j=1}^{l} \Lambda _{i_j}$ $( l \in \Zint_{\geq 0}, \; i_j \in \{ 1,\dots ,N-1\}$ $(j=1,\dots ,l))$

We repeatedly apply Sublemma \ref{sublemma1} for $J_{\lambda }\tilde{e}_{\nu }$.
Then $J_{\lambda }\tilde{e}_{\nu }$ is expressed as the linear combination of $J_{\nu '}$, where $\nu' = \lambda +\sum_{j=1}^{l} w_j(\Lambda _{i_j})$ for some $w_j \in W$ $(j=1,\dots l)$.
From Sublemma \ref{sublemma2}, we have $|\nu' - \lambda| \leq |\sum_{j=1}^{l} \Lambda _{i_j}|=|\nu |$. Applying Sublemma \ref{sublemma4} for $\mu $ and $\nu$, we obtain Lemma \ref{anallem2}.
\end{proof}

\begin{prop} \label{analprop1}
Let $|p|<1$ and $\lambda \in P_+$. Write $(\sum_{k=1}^{\infty} T^{(k)}p^{k})\Tilde{J}_{\lambda}=\sum_{\mu \in P_+, \lambda -\mu \in Q}\Tilde{t}_{\lambda,\mu}\Tilde{J}_{\mu}$, where $\Tilde{J}_{\lambda}$ is the normalized $A_{N-1}$-Jacobi polynomial. For each $C$ such that $C>1$ and $C|p|<1$, there exists a number $C''\in \Rea _{>0}$ such that $|\Tilde{t}_{\lambda,\mu}| \leq C''(C|p|)^{(|\lambda -\mu|+1)/2\sqrt{2}}$ for all $\mu \in P_+$. 

\end{prop}
\begin{proof}
Since the normalized Jacobi polynomials form the complete orthonormal system with respect to the inner product $(,)_{\Delta }$, we have $\Tilde{t}_{\lambda,\mu}= ( (\sum_{k=1}^{\infty} T^{(k)}p^{k})\Tilde{J}_{\lambda}, \Tilde{J}_{\mu} ) _{\Delta }$.

We fix $\lambda , \mu \in P_+$.
Let $m$ be the smallest integer which is greater or equal to $|\lambda -\mu|/\sqrt{2}$.
If $k<m$, then we have $ ( T^{(k)}p^{k}\Tilde{J}_{\lambda} ,  \Tilde{J}_{\mu} ) _{\Delta }=0$ by Lemmas \ref{anallem1}, \ref{anallem2} and the orthogonality.
Therefore we have
\begin{align*}
& |\Tilde{t}_{\lambda,\mu}|= \left|(\Tilde{J}_{\mu}, (\sum_{k=1}^{\infty} T^{(k)}p^{k})\Tilde{J}_{\lambda} ) _{\Delta } \right|\\
& = \left|( \Tilde{J}_{\mu},(\sum_{k=m}^{\infty} T^{(k)}p^{k})\Tilde{J}_{\lambda} ) _{\Delta } \right|\\
& = \left|\int _T \sum_{k=m}^{\infty}\sum_{\alpha \in \Delta_{+}}t_k(\langle \alpha, h\rangle )p^{k}\Tilde{J}_{\lambda} \overline{\Tilde{J}_{\mu}} \Delta ^2 d\mu \right|\\
& \leq \sup_{h\in T} \left| \sum_{k=m}^{\infty}\sum_{\alpha \in \Delta_{+}}t_k(\langle \alpha, h\rangle )p^{k} \right| \int_T |\Tilde{J}_{\lambda} \Tilde{J}_{\mu} \Delta ^2 | d\mu \\
& \leq \frac{k_{\alpha}( k_{\alpha}-1)N(N-1)}{2}\sum_{k \geq m}t_k|p|^{k} \cdot \sqrt{ \left(\int_T |\Tilde{J}_{\lambda}  \Delta |^2 d\mu \right)\left(\int_T |\Tilde{J}_{\mu}  \Delta |^2 d\mu \right)}\\
& \leq \frac{k_{\alpha}( k_{\alpha}-1)N(N-1)}{2}\sum_{k \geq m}t_k|p|^{k}.
\end{align*}

Similarly, we have $|\Tilde{t}_{\lambda,\lambda}| \leq \frac{k_{\alpha}( k_{\alpha}-1)N(N-1)}{2}\sum_{k \geq 1}t_k|p|^{k}$.

Since the convergence radius of the series $\sum_n t_n p^n$ is equal to $1$, we obtain that there exists a number $C''\in \Rea _{>0}$ such that $|\Tilde{t}_{\lambda,\lambda}| \leq C''(C|p|)$ and $|\Tilde{t}_{\lambda,\mu}| \leq C''(C|p|)^{|\lambda -\mu|/\sqrt{2}}$ for $\lambda \neq \mu$.
Hence we have the Proposition.
\end{proof}

\begin{prop} \label{analprop11}
Let $D$ be a positive number. We suppose $\dist(\zeta, \sigma (\Tilde{H}_0))\geq D$. Write $(\Tilde{T}(p)-\zeta )^{-1}  \Tilde{J}_{\lambda}=\sum_{\mu}t_{\lambda,\mu}\Tilde{J}_{\mu}$, where $(\Tilde{T}(p)-\zeta )^{-1}$ is defined in (\ref{Neus}).
For each $\lambda \in P_+$ and $C\in \Rea _{>1}$, there exists $C'\in \Rea _{>0}$ and $p_0\in \Rea _{>0}$ which do not depend on $\zeta $ (but depend on $D$) such that $t_{\lambda,\mu}$ satisfies 
\begin{equation}
|t_{\lambda,\mu}| \leq C'(C|p|)^{\frac{|\lambda -\mu|}{2N\sqrt{2}}},
\label{proptlm}
\end{equation}
 for all $p,$ $\mu$ s.t. $|p|<p_0$ and $\mu \in P_+$.
\end{prop}
\begin{proof}
Let us recall that the operator $(\Tilde{T}(p)-\zeta )^{-1}$ is defined by the Neumann series (\ref{Neus}).

We fix the number $D(\in \Rea_{>0})$ and set $X:=(\zeta -\Tilde{H}_0)^{-1}(\sum_{k=1}^{\infty} T^{(k)}p^{k})$.
From the expansion (\ref{Neus}) and Proposition \ref{tkineq}, there exists a number $p_1\in \Rea_{>0}$ such that the inequality $\| X \| <1/2$ holds for $p$ ($|p|<p_1$) and $\zeta$ ($\dist(\zeta, \sigma (\Tilde{H}_0))>D)$. In this case, we have $(\sum_{i=0}^{\infty}X^i)(\Tilde{H}_0-\zeta )^{-1}= (\Tilde{T}(p)-\zeta )^{-1}$.
We write $X\Tilde{J}_{\lambda}=\sum_{\mu}\check{t}_{\lambda,\mu}\Tilde{J}_{\mu}$.

For the series $\sum_{\mu} c_{\mu}\Tilde{J}_{\mu}$, write $\sum_{\mu} c'_{\mu}\Tilde{J}_{\mu}=(\Tilde{H}_0-\zeta )^{-1} \sum_{\mu}c_{\mu}\Tilde{J}_{\mu}$.
We have $|c'_{\mu}| \leq D^{-1}|c_{\mu}|$ for each $\mu$.
Combining with Proposition \ref{analprop1}, we obtain that for each $C$ such that $C>1$ and $\frac{C+1}{2} p_1<1$. there exists $C''\in \Rea _{>0}$ which does not depend on $\zeta $ but $D$ such that $|\check{t}_{\lambda,\mu}| \leq C''(\frac{C+1}{2}|p|)^{(|\lambda -\mu|+1)/2\sqrt{2}}$ if $|p|<p_1$.

To obtain Proposition \ref{analprop11}, we use the method of majorants.

We introduce the symbol $\bold{e}_{\lambda }$ $(\lambda \in (\frac{\Zint}{N} )^N)$ to avoid inaccuracies. We remark that $P_+ \subsetneq  (\frac{\Zint}{N} )^N$.
We will apply the method of majorants for the formal series  $\sum_{\mu\in (\frac{\Zint}{N} ) ^N} c_{\mu}\bold{e}_{\mu}$ instead of $\sum_{\mu\in P_+} c_{\mu}\bold{e}_{\mu}$ .

For the formal series $\sum_{\mu\in (\frac{\Zint}{N} ) ^N} c_{\mu}\bold{e}_{\mu}$, we define the partial ordering $\Tilde{\leq }$ by the rule
$$
\sum_{\mu\in (\frac{\Zint}{N} ) ^N} c^{(1)}_{\mu}\bold{e}_{\mu} \Tilde{\leq } \sum_{\mu\in (\frac{\Zint}{N} ) ^N} c^{(2)}_{\mu}\bold{e}_{\mu}\; 
\Leftrightarrow\; \forall \mu, \;  |c^{(1)}_{\mu}|\leq |c^{(2)}_{\mu}| .
$$

We will later consider the case that each $c_{\mu}^{(i)}$ $(\mu \in (\frac{\Zint}{N} )^N, \; i=1,2)$ is expressed as the infinite sum.
If one shows the absolute convergence of $c_{\mu}^{(2)}$ for each $\mu$, one has the absolute convergence of $c_{\mu}^{(1)}$ for each $\mu$ by the majorant.

We set $X\bold{e}_{\lambda}=\sum_{\mu}\check{t}_{\lambda,\mu}\bold{e}_{\mu}$, where the coefficients $\check{t}_{\lambda,\mu}$ were defined by $X\Tilde{J}_{\lambda}=\sum_{\mu}\check{t}_{\lambda,\mu}\Tilde{J}_{\mu}$.

Our goal is to show (\ref{proptlm}) for $t_{\lambda,\mu}$ s.t. $\sum_{\mu} t_{\lambda,\mu} \Tilde{J}_{\mu }=  (\sum_{i=0}^{\infty} X^i )(\Tilde{H}_0-\zeta )^{-1}\Tilde{J}_{\lambda}$.
Since $(\Tilde{H}_0-\zeta )^{-1}\Tilde{J}_{\lambda}= (E_{\lambda} -\zeta )^{-1}\Tilde{J}_{\lambda}$ and $|(E_{\lambda} -\zeta )^{-1}|\leq D^{-1}$, it is enough to show that 
there exists $C^{\star} \in \Rea _{>0}$ and $p_0 \in \Rea _{>0}$ which do not depend on $\zeta $ (but depend on $D$) such that $t^{\star}_{\lambda,\mu}$ is well-defined by  $\sum_{k=0}^{\infty} X^k \bold{e}_{\lambda}= t^{\star}_{\lambda,\mu}\bold{e}_{\mu}$ and satisfies 
\begin{equation}
|t^{\star}_{\lambda,\mu}| \leq C^{\star}\left(C|p|\right)^{\frac{|\lambda -\mu|}{2N\sqrt{2}}},
\label{proptlms}
\end{equation}
 for all $p,$ $\mu$ s.t. $|p|<p_0$ and $\mu \in P_+$.

We set
$$
Y \bold{e}_{\lambda} := \sum_{\mu, \lambda -\mu \in \Zint ^N}y_{\lambda, \mu} \bold{e}_{\mu}:=\sum_{\mu,\lambda -\mu \in \Zint ^N}C''\left(\frac{C+1}{2}p\right)^{(|\mu -\lambda |+1)/2\sqrt{2}} \bold{e}_{\mu}.
$$
$$
Z \bold{e}_{\lambda} := \sum_{\mu,\lambda -\mu \in \Zint ^N}z_{\lambda, \mu} \bold{e}_{\mu}:=\sum_{\mu,\lambda -\mu \in \Zint ^N}C''\left(\frac{C+1}{2}p\right)^{\sum_{i=1}^{N}\frac{(|\mu_i- \lambda_i |+1)}{2N\sqrt{2}}} \bold{e}_{\mu}.
$$
We have the inequality
$$
X \bold{e}_{\lambda} \Tilde{\leq } Y \bold{e}_{\lambda} \Tilde{\leq } Z \bold{e}_{\lambda}.
$$
Let $k \in \Zint_{\geq 1}$. If the coefficients of $Z^k \bold{e}_{\lambda}$ w.r.t the basis $\{ \bold{e}_{\mu} \}$ converge absolutely, then the coefficients of the series $X^k \bold{e}_{\lambda}$ and $Y^k \bold{e}_{\lambda}$ are well-defined and we have
$$
X^k \bold{e}_{\lambda} \Tilde{\leq } Y^k \bold{e}_{\lambda} \Tilde{\leq } Z^k \bold{e}_{\lambda}.
$$
From the equality $Z^k \bold{e}_{\lambda}= \sum_{\nu^{(1)}, \dots, \nu^{(k-1)}} z_{\lambda , \nu^{(1)}}z_{\nu^{(1)},\nu^{(2)}}\dots z_{\nu^{(k-1)},\mu} \bold{e}_{\mu}$ and the property $z_{\lambda ,\mu }=z_{0,\mu -\lambda}$, we have
\begin{align*}
& Z^k \bold{e}_{\lambda}
=  \sum_{\mu, \lambda -\mu \in \Zint ^N}\frac{1}{(2\pi\sqrt{-1})^N}\oint _{|s_1 |=1}\dots \oint _{|s_N |=1} \\
& \left(\sum_{\nu =(\nu_1, \dots ,\nu_N) \in  \Zint ^N}z_{0,\nu }s_1^{\nu_1}\dots s_N^{\nu_N}\right)^k s_1^{\lambda_1 -\mu _1-1}\dots  s_N^{\lambda_N -\mu _N-1}ds_1 \dots ds_N  \bold{e}_{\mu}\\
&=\sum_{\mu, \lambda -\mu \in \Zint ^N}\frac{1}{(2\pi\sqrt{-1})^N} \prod_{i=1}^N \oint_{|s_i|=1}\left(\sum_{\nu \in \Zint} (C'')^{\frac{1}{N}}\left(\frac{C+1}{2}p\right)^{ \frac{(|\nu |+1)}{2N\sqrt{2}}}s_i^{\nu }\right)^k s_i^{\lambda_i-\mu_i-1}ds_i \bold{e}_{\mu}.
\end{align*}
Set $\Tilde{p}= (\frac{C+1}{2}p)^{\frac{1}{2N\sqrt{2}}}$, we have
\begin{align*}
& \sum_{k=1}^{\infty}Z^k \bold{e}_{\lambda} \\
& = \sum_{k=1}^{\infty}\sum_{\mu, \lambda -\mu \in \Zint ^N}\frac{1}{(2\pi\sqrt{-1})^N} \prod_{i=1}^N \oint_{|s_i|=1}\left(\sum_{\nu \in \Zint} (C'')^{\frac{1}{N}}\Tilde{p}^{(|\nu |+1)}s_i^{\nu }\right)^k s_i^{\lambda_i-\mu_i-1}ds_i \bold{e}_{\mu} \\
& \Tilde{\leq }\sum_{\mu, \lambda -\mu \in \Zint ^N}\frac{1}{(2\pi\sqrt{-1})^N} \prod_{i=1}^N \oint_{|s_i|=1}\sum_{k=1}^{\infty}\left(\sum_{\nu \in \Zint} (C'')^{\frac{1}{N}}\Tilde{p}^{(|\nu |+1)}s_i^{\nu }\right)^k s_i^{\lambda_i-\mu_i-1}ds_i \bold{e}_{\mu} \\
& = \sum_{\mu, \lambda -\mu \in \Zint ^N} Z_{\lambda ,\mu} \bold{e}_{\mu}, 
\end{align*}
where 
\begin{align}
& Z_{\lambda ,\mu}= \frac{1}{(2\pi\sqrt{-1})^N} \prod_{i=1}^N \oint_{|s_i|=1}\frac{(C'')^{\frac{1}{N}}(\Tilde{p}-\Tilde{p}^3)s_i^{\lambda_i-\mu_i-1}ds_i}{(1-\Tilde{p}s_i)(1-\Tilde{p}s_i^{-1})-(C'')^{\frac{1}{N}}(\Tilde{p}-\Tilde{p}^3)} .
\label{zlm}
\end{align}
Remark that we used the inequality $\sum_{k=1}^{\infty} \prod_{i=1}^{N} (a_i)^k \leq \prod_{i=1}^{N}\sum_{k=1}^{\infty} (a_i)^k$ for $0<a_1,\dots ,a_N<1$ and the formula $\sum_{n\in \Zint}q^{|n|+1}x^n=\frac{q-q^3}{(1-qx)(1-qx^{-1})}$. The equality (\ref{zlm}) makes sense for $\tilde{p}<p_2$, where $p_2$ is the positive number satisfying the inequalities $p_2<1$, $\frac{(C'')^{\frac{1}{N}}|p_2-p_2^3|}{(1-p_2)^2}<1$ and $(C'')^{\frac{1}{N}}p_2<1$.

Therefore each coefficient of $\sum_{k=1}^{\infty} Z^k \bold{e}_{\lambda}$ w.r.t the basis $\{ \bold{e}_{\mu} \}$ converges absolutely. Hence the following inequality makes sense,
$$
\sum_{k=0}^{\infty} X^k \bold{e}_{\lambda}  \Tilde{\leq }  \bold{e}_{\lambda}+ \sum_{k=1}^{\infty} Z^k \bold{e}_{\lambda} \Tilde{\leq } \bold{e}_{\lambda}+  \sum_{\mu, \lambda-\mu \in \Zint^{N}} Z_{\lambda ,\mu} \bold{e}_{\mu}.
$$

Let $s(\Tilde{p})$ be the solution of the equation $(1-\Tilde{p}s)(1-\Tilde{p}s^{-1})-(C'')^{\frac{1}{N}}(\Tilde{p}-\Tilde{p}^3)=0$ on $s$ satisfying $|s(\Tilde{p})|<1$. Then $s(\Tilde{p})$ is holomorphic in $\Tilde{p}$ near $0$ and admits the expansion $s(\Tilde{p})= \Tilde{p}+ c_2\Tilde{p}^2+\dots $.
We have
\begin{equation}
\frac{1}{(2\pi\sqrt{-1})}\oint_{|s|=1}\frac{(C'')^{\frac{1}{N}}(\Tilde{p}-\Tilde{p}^3)s^{n-1}ds}{(1-\Tilde{p}s)(1-\Tilde{p}s^{-1})-(C'')^{\frac{1}{N}}(\Tilde{p}-\Tilde{p}^3)}= \Tilde{p}f(\Tilde{p})s(\Tilde{p})^{|n|},
\label{integTp}
\end{equation}
where $f(\Tilde{p})$ is a holomorphic function defined near $\Tilde{p}=0$.
For the $n\geq 0$ case, we have the relation (\ref{integTp}) by calculating the residue around $s=s(\Tilde{p})$. For the $n<0$ case, we need to change the variable $s \rightarrow s^{-1}$ and calculate the residue around $s=s(\Tilde{p})$.

The coefficient of $\bold{e}_{\mu}$ on the series $\sum_{k=0}^{\infty} X^k\bold{e}_{\lambda}$ satisfying $\lambda -\mu \not \in Q$ has to be zero from the definition of $X$.
By the inequality $|s(\Tilde{p})|^{|\nu_1|+\dots +|\nu_N|}\leq |s(\Tilde{p})|^{\sqrt{\nu_1^2+\dots +\nu_N^2}}$, we have
$$
\sum_{k=0}^{\infty}X^k \bold{e}_{\lambda} \Tilde{\leq }\bold{e}_{\lambda} + \sum_{\mu, \lambda -\mu\in Q} \left( \prod_{i=1}^N (\Tilde{p}f(\Tilde{p})) s(\Tilde{p})^{|\lambda_i -\mu _i |} \right) \bold{e}_{\mu} \Tilde{\leq }\bold{e}_{\lambda} + \sum_{\mu, \lambda -\mu\in Q} (\Tilde{p}f(\Tilde{p}))^Ns(\Tilde{p})^{|\lambda -\mu |} \bold{e}_{\mu}.
$$
for $|\Tilde{p}|<p_3$ and $p_3$: a sufficiently small positive number.

Combining with the relation $\Tilde{p}=\left(\frac{C+1}{2}p\right)^{\frac{1}{2N\sqrt{2}}}$, the inequality $\frac{C+1}{2}<C$, and the expansion $s(\Tilde{p})= \Tilde{p}+ c_2\Tilde{p}^2+\dots $, we obtain (\ref{proptlms}) and the proposition.
\end{proof}

\begin{prop} \label{analprop2}
Let $a_i \in \sigma (\Tilde{H}_0)$ and $\Gamma _i$ be a circle in $\Cplx$ which contains only one element $a_i$ of $\sigma (\Tilde{H}_0)$ inside it.
Let $\lambda \in P_+$ satisfying $\Tilde{H}_0\Tilde{J}_{\lambda}= a_i \Tilde{J}_{\lambda}$.

We set $P_i(p)= -\frac{1}{2\pi \sqrt{-1}}\int_{\Gamma_i}(\Tilde{T}(p)-\zeta )^{-1} d\zeta $ and write $P_i(p) \Tilde{J}_{\lambda}=\sum_{\mu}s_{\lambda,\mu}\Tilde{J}_{\mu}$. 

For each $C\in \Rea _{>1}$, There exists $C'\in \Rea _{>0}$ and $p_{\ast }\in \Rea _{>0}$ such that $s_{\lambda,\mu}$ satisfies 
\begin{equation}
|s_{\lambda,\mu}| \leq C'(C|p|)^{\frac{|\lambda -\mu|}{2N\sqrt{2}}},
\label{propslm}
\end{equation}
 for all $p,$ $\mu$ s.t. $|p|<p_{\ast }$ and $\mu \in P_+$.
\end{prop}
\begin{proof}

Since the spectrum $\sigma (\Tilde{H}_0)$ is discrete,
there exists a positive number $D$ such that $\inf_{\zeta \in \Gamma_i}\dist(\zeta, \sigma (\Tilde{H}_0)) \geq D$.

We write $(\Tilde{T}(p)-\zeta )^{-1} \Tilde{J}_{\lambda}=\sum_{\mu}t_{\lambda,\mu}(\zeta )\Tilde{J}_{\mu}$.
From Proposition \ref{analprop11}, we obtain that for each $C\in \Rea _{>1}$, there exists $C_{\ast }\in \Rea _{>0}$ and $p_{\ast }\in \Rea _{>0}$ which does not depend on $\zeta (\in \Gamma_i) $  such that $|t_{\lambda,\mu}(\zeta )| \leq C_{\ast }(C|p|)^{\frac{|\lambda -\mu|}{2N\sqrt{2}}}$ for all $p,$ $\mu$ s.t. $|p|<p_{\ast }$ and $\mu \in P_+$.

Let $L$ be the length of the circle $\Gamma_i$ and write $-\frac{1}{2\pi \sqrt{-1}} \int_{\Gamma_i}(\Tilde{T}(p)-\zeta )^{-1} d\zeta  \Tilde{J}_{\lambda}=\sum_{\mu}s_{\lambda,\mu}\Tilde{J}_{\mu}$.
By integrating $\sum_{\mu}t_{\lambda,\mu}(\zeta )\Tilde{J}_{\mu}$ over the circle $ \Gamma_i$, we have $|s_{\lambda,\mu}| \leq \frac{L}{2\pi} C_{\ast }(C|p|)^{\frac{|\lambda -\mu|}{2N\sqrt{2}}}$ for all $p,$ $\mu$ s.t. $|p|<p_{\ast }$ and $\mu \in P_+$.

Therefore we have Proposition \ref{analprop2}.
\end{proof}

\begin{prop}  \label{analprop3}
Let $\mu \in P_+$ and $c_{\mu}$ be a number satisfying $|c_{\mu}|<a|p|^{b|\mu|}$ $(|\mu|>M)$ for some $a,b>0$ and $M \in \Zint $.
The function $\sum_{\mu}c_{\mu}\Tilde{J}_{\mu}$ is holomorphic  when $|Im x_j |$ $(j=1,\dots ,N)$ and $|p|$ are sufficiently small.
\end{prop}
\begin{proof}
Since $z_i=e^{\sqrt{-1}x_i}$, it is enough to show that the function $\sum_{\mu}c_{\mu}\Tilde{J}_{\mu}$ is holomorphic  when $|p|$ is sufficiently small and $1/2<|z_j|<2$ $(j=1,\dots ,N)$.

We count roughly the number of the elements of $P_+$ of a given length. The rough estimate is given by $\# \{ \lambda \in P_+ \: | \: (\lambda |\lambda )=l\} \leq (2lN)^N $. We will use this in the inequality (\ref{ineq}).

In the proof, we will use the notations and the results written in section \ref{Jack}. In  section \ref{Jack}, there are parameters $r$ and $C_0$. We fix $r=2$ and $C_0=2$. There is another number $A$ defined in  section \ref{Jack}.

We have
\begin{align}
& \left|\sum_{\mu \in P_+, |\mu |\geq M\atop{1/2\leq |z_i|\leq 2}} c_{\mu}\Tilde{J}_{\mu}(z_1,\dots ,z_N)\right| \leq \sum_{\mu \in P_+, |\mu |\geq M\atop{1/2\leq |z_i|\leq 2}}a|p|^{b|\mu|}|\Tilde{J}_{\mu}(z_1,\dots ,z_N)| \label{ineq} \\
& \leq \sum_{\mu \in P_+, |\mu |\geq M}
a2^{(N-1)\sqrt{2({\mu }|{\mu })}}2^{\sqrt{N(\mu | \mu )}} |p|^{b\sqrt{(\mu | \mu )}} \notag \\
& \leq \sum_{n\geq M, n\in \Zint/N}aA(2nN)^N(2^{(N-1)\sqrt{2}+\sqrt{N}}|p|^b)^{\sqrt{n}}.
\notag
\end{align}
by the formulae (\ref{lineq}), (\ref{Jineq}).

If $2^{(N-1)\sqrt{2}+\sqrt{N}}|p|^b<1$ then the bottom part of the inequality converges. We choose a positive number $p_0$ which satisfies $2^{(N-1)\sqrt{2}+\sqrt{N}}p_0^b<1$. Then the series $|\sum_{\mu \in P_+, |\mu |\geq M} c_{\mu}\Tilde{J}_{\mu}(z_1,\dots ,z_N)|$ is uniformly bounded and uniformly absolutely converges for $|p|<p_0$ and $1/2\leq|z_i| \leq 2$ $(i=1,\dots ,N)$.
Since the functions $\Tilde{J}_{\mu}(z_1,\dots ,z_N)$ are holomorphic, we have the holomorphy of the function $\sum_{\mu}c_{\mu}\Tilde{J}_{\mu}(z_1,\dots ,z_N)$ by the Weierstrass's theorem.
\end{proof}

Combining Propositions \ref{analprop2} and \ref{analprop3}, we have Proposition \ref{analthm}.

From Propositions \ref{analprop11} and \ref{analprop3}, the function $(\tilde{T}(p)-\zeta )^{-1}\Tilde{J}_{\lambda }$ is real-holomorphic on $(x_1,\dots ,x_N)$ if $|p|$ is sufficiently small. From Proposition \ref{prop:actwelldef}, the operators $H^{(j)}(p)$ $(j=1,\dots ,N)$ act well-definedly on the function $(\tilde{T}(p)-\zeta)^{-1}\Tilde{J}_{\lambda }$ and we have $H^{(j)}(p)  (\tilde{T}(p)-\zeta)^{-1}\Tilde{J}_{\lambda } \in C^{\omega}(T)^W$. It follows from the commutativity of the operators $\tilde{T}(p)$ and $H^{(j)}(p)$ (\ref{THcomm}) that
\begin{align*}
& H^{(j)}(p)\Tilde{J}_{\lambda }= H^{(j)}(p)(\tilde{T}(p)-\zeta)(\tilde{T}(p)-\zeta)^{-1}\Tilde{J}_{\lambda } \\
& = (\tilde{T}(p)-\zeta)H^{(j)}(p)(\tilde{T}(p)-\zeta)^{-1}\Tilde{J}_{\lambda }.
\end{align*}
Hence we have $H^{(j)}(p)(\tilde{T}(p)-\zeta)^{-1}\Tilde{J}_{\lambda }=(\tilde{T}(p)-\zeta)^{-1}H^{(j)}(p)\Tilde{J}_{\lambda }$.

By integrating it on the variable $\zeta$ over the circle $\Gamma _i$, we have $H^{(j)}(p)P_i(p)\Tilde{J}_{\lambda }=P_i(p) H^{(j)}(p)\Tilde{J}_{\lambda }$.

Therefore we have Proposition \ref{prop:commHP}.

\appendix
\section{Jack polynomial and special functions}
\subsection{Jack polynomial and $A_{N-1}$-Jacobi polynomial}\label{Jack}

We will see the relationship between the Jack polynomial and the $A_{N-1}$-Jacobi polynomial.

Let $\mathcal{M}_N$ be the set of partitions with at most $N$ parts, \emph{i.e.}, 
$\mathcal{M}_N:= \{ \underline{\lambda } =(\underline{\lambda } _1, \dots ,\underline{\lambda } _N ) \: | \:  \underline{\lambda }_i -\underline{\lambda }_{i+1} \in \Zint _{\geq 0}, ( i= 1,\dots ,N-1), \: \underline{\lambda }_N \in \Zint_{\geq 0}\}$.
We set $\mathcal{M}^0_N:=  \{ \underline{\lambda } =(\underline{\lambda } _1, \dots ,\underline{\lambda } _N ) \in \mathcal{M}_N  \: | \: \underline{\lambda }_N=0 \}$.
The Jack polynomial $J^{\star}_{\underline{\lambda } }(z_1, \dots , z_N)$ $(\underline{\lambda } \in \mathcal{M}_N)$ is a symmetric polynomial of variables $(z_1, \dots , z_N)$ which is a eigenfunction of the gauge-transformed Hamiltonian $\Tilde{H}_0$ (\ref{gtHam}).

Let $m^{\star}_{\underline{\lambda } }$ be the monomial symmetric polynomial.
The Jack polynomial admits the following expansion,
\begin{equation}
J^{\star}_{\underline{\lambda }}= m^{\star}_{\underline{\lambda }}+ \sum_{\underline{\mu} \prec \underline{\lambda }}u_{\underline{\lambda } \underline{\mu}}m^{\star}_{\underline{\mu}},
\end{equation}
where the dominant ordering of $\mathcal{M}_N$ is given by $\underline{\lambda }  \preceq \underline{\mu} \Leftrightarrow \sum_{j=1}^{i} \underline{\lambda }_j \leq \sum_{j=1}^{i} \underline{\mu}_j \: (i=1,\dots ,N-1), \sum_{j=1}^{N} \underline{\lambda }_j = \sum_{j=1}^{N} \underline{\mu}_j$.

We see the correspondence between the Jack polynomial and the $A_{N-1}$-Jacobi polynomial.
Let $J^{\star}_{\underline{\lambda } }(z_1, \dots , z_N)$ $(\underline{\lambda } \in \mathcal{M}_N)$ be a Jack polynomial. We set $|\underline{\lambda } |_{\star}=\sum_{i=1}^{N} \underline{\lambda } _i$ and $\lambda= \sum_{i=1}^{N} (\underline{\lambda }_i -|\underline{\lambda } |_{\star}/N)\epsilon _i$. Then $\lambda \in P_+$, where $P_+$ is the set of dominant weights of type $A_{N-1}$.
The function $(z_1 \dots  z_N)^{ -|\underline{\lambda } |_{\star}/N}J^{\star}_{\underline{\lambda } }(z_1, \dots , z_N)$ is precisely the $A_{N-1}$-Jacobi polynomial $J_{\lambda}$.

By this correspondence, the Jack polynomial $J^{\star}_{\underline{\lambda } }(z_1, \dots , z_N)$ $(\underline{\lambda } \in \mathcal{M}^0 _N)$ corresponds with the $A_{N-1}$-Jacobi polynomial $J_{\lambda}$ $(\lambda \in P_+)$ one-to-one.

Let $\underline{\lambda } $ be an element in $\mathcal{M}^0 _N$ and $\lambda $ be the corresponding element in $P_+$.
Since $(\lambda | \lambda) \geq (\underline{\lambda }_1 -|\underline{\lambda } |_{\star}/N)^2+(|\underline{\lambda } |_{\star}/N)^2\geq (\underline{\lambda } _1)^2/2 \geq \frac{|\underline{\lambda } |_{\star}^2}{2(N-1)^2}$, we have
\begin{equation}
|\underline{\lambda } |_{\star}\leq (N-1)\sqrt{2(\lambda | \lambda )}.
\label{lineq}
\end{equation}

Let us remind the Cauchy formula for the Jack polynomial.
\begin{equation}
 \prod_{1\leq i,j\leq N}(1-\kappa X_iY_j)^{-\beta}
= \sum_{\underline{\lambda } \in \mathcal{M}_N} \kappa^{|\underline{\lambda }|_{\star}}J^{\star}_{\underline{\lambda }}(X)J^{\star}_{\underline{\lambda }}(Y)j_{\underline{\lambda }}^{-1},
\label{Cauchy}
\end{equation}
where 
\begin{equation}
  0\leq j_{\underline{\lambda }}=\prod_{s\in \underline{\lambda }}
\frac{ a(s)+\beta l(s)+1}{a(s)+\beta l(s)+\beta}\leq 1,
\end{equation}
due to $\beta\geq 1$.
$a(s)$ is the arm-length and $l(s)$ is the leg-length.

Since p.379 of Macdonald's book (\cite{Mac}), we have
$J^{\star}_{\underline{\lambda }}= m^{\star}_{\underline{\lambda }}+ \sum_{\underline{\mu} \prec \underline{\lambda }}u_{\underline{\lambda } \underline{\mu}}m^{\star}_{\underline{\mu}}$ with $u_{\underline{\lambda } \underline{\mu}}>0$ if $\beta >0$.
Hence we have
\begin{equation}
J_\lambda= m_{\lambda}+ \sum_{\lambda -\mu \in Q_+}u_{\lambda \mu}m_{\mu},
\end{equation}
with $u_{\lambda \mu}>0$.
Let $r$ be a real number greater than $1$.

If $1/r <|z_i|<r$ for all $i$ then $|m_{\lambda}(z_1, \dots ,z_N)|\leq r^{\sum_{i=1}^N {|(\lambda| \epsilon_i)|}}m_{\lambda}(1) \leq r^{\sqrt{N(\lambda| \lambda)}} m_{\lambda}(1)$.
Therefore we have
\begin{equation}
 0\leq |J_\lambda(z)|\leq  r^{\sqrt{N(\lambda| \lambda)}}J_\lambda(1)
\end{equation}
on $1/r <|z_i|<r$ for all $i$.

By setting $X_i=Y_j=1$ in (\ref{Cauchy}), we have
\begin{equation}
 (1-\kappa)^{-\beta N^2}
=\sum_{n\in \Zint_{\geq 0}}\kappa^{n}c_n
=\sum_{\underline{\lambda } \in \mathcal{M}_N}\kappa^{|\underline{\lambda }|_{\star}}J^{\star}_{\underline{\lambda }}(1)^2j_{\underline{\lambda }}^{-1},
\end{equation}
where $c_n=\frac{\Gamma(\beta N^2+n+1)}{\Gamma(\beta N^2+1)\Gamma(n+1)}$.
For each $\beta \geq 1$ and $C_0 >1$, there exists a positive number $A$ such that $c_n<A^2C_0^{2n}$ for all $n \in \Zint _{\geq 0}$.
Thus
\begin{equation}
J^{\star}_{\underline{\lambda }}(1)^2 \leq \sum_{|\underline{\lambda }|_{\star}=n}J^{\star}_{\underline{\lambda }}(1)^2 j_{\underline{\lambda }}^{-1}<A^2 C_0^{2n}.
\end{equation}
By the inequality (\ref{lineq}), we have
\begin{equation}
|J_{\lambda }(1)| <A C_0^{(N-1)\sqrt{2({\lambda}|{\lambda})}}.
\end{equation}
The square of the norm of $J^{\star}_{\underline{\lambda }}$ is
\begin{equation}
\|J^{\star}_{\underline{\lambda }}\|^2=\prod_{i<j}
\frac
{\Gamma (\xi_i-\xi_j+\beta)\Gamma (\xi_i-\xi_j-\beta+1)}
{\Gamma (\xi_i-\xi_j)\Gamma (\xi_i-\xi_j+1)}
\end{equation}
, where $\xi_i=\underline{\lambda }_i +\beta(N-i)$. (See (\cite{Mac} p.383))
If $\beta \geq 1$ then we have $\|J^{\star}_{\underline{\lambda }}\|^2 \geq 1$ because of the convexity of the function $\log \Gamma (x)$.
Therefore we have $\|J_\lambda\|^2 \geq  1$.
Generally we have for $r>1$
\begin{equation}
\max_{1/r \leq |z_i| \leq r} |\Tilde{J}_\lambda(z)| 
\leq 
\max_{1/r \leq |z_i| \leq r} |J_\lambda(z)| 
\leq 
AC_0^{(N-1)\sqrt{2({\lambda}|{\lambda})}}r^{\sqrt{N(\lambda| \lambda)}},
\label{Jineq}
\end{equation}
where $\Tilde{J}_\lambda(z)$ is the normalized $A_{N-1}$-Jacobi polynomial.

\subsection{Special functions} \label{fns}

We define some functions needed in this article.
\begin{equation}
\theta_1(x):=2\sum_{n=1}^{\infty} (-1)^{n-1} \exp (\tau \pi \sqrt{-1}(n-1/2)^2) \sin(2n-1)\pi x, \label{th1}
\end{equation}
\[
\theta (x):=\frac{\theta_1(x)}{\theta_1 '(0)}
\]
\[
\wp (x; \omega _1, \omega _3):=\frac{1}{z^2}+\sum_{(m,n) \setminus \{ (0,0)\} \in \Zint ^2} 
\left(
\frac{1}{(z+2m\omega _1+2n\omega _3)^2}- \frac{1}{(2m\omega _1+2n\omega _3)^2}
\right)
\]
\[
\wp (x):=\wp (x; \pi, \pi \tau ).
\]
We have
\begin{equation}
\wp (x)=
\frac{1}{4\sin^2 (x/2)}- \frac{1}{12} -2 \sum_{n=1}^{\infty} \frac{np ^{n}}{1-p ^{n}} (\cos n x -1), \label{wpth}
\end{equation}
where $p= \exp(2\tau \pi \sqrt{-1})$.

\vspace{.15in}
{\bf Acknowledgment}
The authors would like to thank Prof.~M. Kashiwara and  Prof.~T. Miwa for discussions and support. Thanks are also due to Dr.~T. Koike and Prof.~T. Oshima.
They thank the referee for valuable comments.
One of the authors (YK) is a Research Fellow of the Japan Society 
for the Promotion of Science.

\end{document}